\documentclass[11pt,a4paper]{article}
\usepackage{epsfig}
\usepackage[]{natbib}
\usepackage{amsmath,amsthm,amssymb,amsfonts}
\usepackage{nccfoots,bbm}
\usepackage{graphicx,color,caption}
\usepackage{nth}
\usepackage[active]{srcltx}
\usepackage{hyperref}      

\newcommand*{\defeq}{\mathrel{\vcenter{\baselineskip0.5ex\lineskiplimit0pt\hbox{\scriptsize.}\hbox{\scriptsize.}}}=}
\newcommand*{\eqdef}{=\mathrel{\vcenter{\baselineskip0.5ex\lineskiplimit0pt\hbox{\scriptsize.}\hbox{\scriptsize.}}}}

\newcommand{\toas}{\stackrel{\rm a.s.}{\longrightarrow}}
\newcommand{\tf}{\mathbf{t}}

\newcommand{\argmin}{\mathop{\mathrm{arg min}}}

\newcommand{\MP}{Mar\v{c}enko-Pastur }
\newcommand{\E}{{\mathbb{E}}}

\newcommand{\Var}{\mbox{\sf Var}}
\newcommand{\Cov}{\mbox{\sf Cov}}
\newcommand{\supp}{{\sf Supp}}
\newcommand{\im}{{\sf Im}}
\newcommand{\re}{{\sf Re}}
\newcommand{\tr}{{\sf Tr}}
\newcommand{\diag}{{\sf Diag}}

\newcommand{\bone}{\mathbbm{1}}
\newcommand{\C}{{\mathbb{C}}}
\newcommand{\R}{{\mathbb{R}}}

\newcommand{\wh}{\widehat}
\newcommand{\wt}{\widetilde}
\newcommand{\be}{\begin{equation}}
\newcommand{\ee}{\end{equation}}
\newtheorem{lemma}{Lemma}[section]

\newtheorem{corollary}{Corollary}[section]

\newtheorem{theorem}{Theorem}[section]

\newtheorem{remark}{Remark}[section]
\numberwithin{equation}{section}
\theoremstyle{plain}

\def\qed{\rule{2mm}{2mm}}
\def\scale{0.65}

\newenvironment{packed_itemize}{
\begin{itemize}
  \setlength{\itemsep}{1pt}
  \setlength{\parskip}{0pt}
  \setlength{\parsep}{0pt}
}{\end{itemize}}

\renewcommand{\baselinestretch}{1.2}

\begin{document}

\title{Spectrum Estimation: A Unified Framework for Covariance Matrix Estimation and
  PCA in Large Dimensions\thanks{Research was partially completed
    while both authors were visiting the Institute for Mathematical Sciences, National University of Singapore in 2012.}}

\bigskip 

\author
 {Olivier Ledoit\\
 Department of Economics \\
 University of Zurich\\
 CH-8032 Zurich, Switzerland\\
\href{mailto:olivier.ledoit@econ.uzh.ch}{olivier.ledoit@econ.uzh.ch}
\and
 Michael Wolf\thanks{
 Research has been supported by 
 the NCCR Finrisk project 
 ``New Methods in Theoretical and Empirical Asset Pricing''.}\\
 Department of Economics \\
 University of Zurich\\
 CH-8032 Zurich, Switzerland\\
 \href{mailto:michael.wolf@econ.uzh.ch}{michael.wolf@econ.uzh.ch}
}

\date
{First version: January 2013\\
This version: July 2013}

\maketitle

\begin{abstract}
Covariance matrix estimation and principal component analysis
(PCA) are two cornerstones of multivariate analysis. Classic
textbook solutions perform poorly when the dimension of the data
is of a magnitude similar to the sample size, or even
larger. In such settings, there is a common remedy for both statistical
problems: {\it nonlinear shrinkage} of the
eigenvalues of the sample covariance matrix. The optimal nonlinear shrinkage formula \mbox{depends} on
unknown population quantities and is thus not available. It is,
however, possible to consistently estimate an oracle nonlinear
shrinkage, which is motivated on asymptotic grounds. A key tool to
this end is consistent estimation of the set of eigenvalues
of the population covariance matrix (also known as the {\it spectrum}), an
interesting and challenging problem in its own right.
Extensive
Monte Carlo simulations demonstrate that our methods have desirable
finite-sample properties and outperform previous proposals.
\end{abstract}

\bigskip

\begin{tabbing}
               \> nonlinear shrinkage, principal component analysis.
\end{tabbing}

\noindent 
JEL CLASSIFICATION NOS: C13.

\newpage 

\section{Introduction}

This paper tackles three important problems in multivariate
statistics: 1) the estimation of the eigenvalues of the covariance
matrix; 2) the estimation of the covariance matrix itself; and
\mbox{3) principal} component analysis (PCA). In many modern
applications, the matrix dimension is not negligible with respect to
the sample size, so textbook solutions based on classic (fixed-dimension)
asymptotics are no longer appropriate. A better-suited framework is {\em
  large-dimensional asymptotics}, where the matrix dimension and the sample
size go to infinity together, while their ratio  --- called the {\em
  concentration} --- converges to a finite, nonzero limit. Under
large-dimensional asymptotics, the sample covariance matrix is no
longer consistent, and neither are its eigenvalues nor its
eigenvectors.

One of the interesting features of large-dimensional asymptotics is
that principal component analysis can no longer be conducted using
covariance matrix eigenvalues. The variation explained by a principal
component is not equal to the corresponding sample eigenvalue and ---
perhaps more surprisingly --- it is not equal to the corresponding
population eigenvalue either. To the best of our knowledge, this fact
has not been noticed before. The variation explained by
a principal component is obtained instead by applying a {\em nonlinear
  shrinkage} formula to the corresponding sample eigenvalue. This
nonlinear shrinkage formula depends on the unobservable population
covariance matrix, but thankfully it can approximated by an {\em
  oracle shrinkage} formula which depends `only' on the unobservable
eigenvalues of the population covariance matrix. This is the
connection with the first of the three problems mentioned above. Once
we have a consistent estimator of the population eigenvalues, we can
use it to derive a consistent estimator of the oracle shrinkage.

The connection with the second problem, the estimation of the whole
covariance matrix, is that the nonlinear shrinkage formula that gives
the variation explained by a principal component also yields the
optimal rotation-equivariant estimator of the covariance matrix
according to the Frobenius norm. Thus, if we can consistently estimate
the population eigenvalues, and if we plug them into the oracle
shrinkage formula, we can address the problems of PCA and covariance
matrix estimation in a unified framework.

It needs to pointed out there that in a rotation-equivariant
framework, consistent (or even improved) estimators of the population
eigenvalues are not available; instead, one needs to retain the sample
eigenvectors. As a consequence, consistent estimation of the
population covariance matrix itself is not possible. Nevertheless, a
rotation-equivariant estimator can still be very useful for practical
applications, as evidenced by the popularity of the previous proposal
of \cite{ledoit:wolf:2004a}. An alternative approach that would allow
for consistent estimation of the population covariance matrix is to
impose additional structure on the estimation problem, such as
sparseness, a graph model, or an (approximate) factor model. But
whether such structure does indeed exist is something that cannot
be verified from the data. Therefore, at least in some applications, a
structure-free approach will be preferred by applied researchers. This
is the problem that we address, aiming to further improve upon
\cite{ledoit:wolf:2004a}. 

Of course estimating population eigenvalues consistently under
large-dimensional asymptotics is no trivial matter. Until recently,
most researchers in the field even feared it might be impossible
because deducing population eigenvalues from sample eigenvalues showed
some   symptoms of {\em ill-posedness}. This means that small
estimation errors in the sample eigenvalues would be amplified by the
specific mathematical structure of the asymptotic relationship between
sample and population eigenvalues. But two recent articles by
\cite{mestre:2008b} and \cite{karoui:2008}  challenged this widely-held
belief and gave some hope that it might be possible after all to
estimate the population eigenvalues consistently. Still, a general
satisfactory solution is not available to date.

The work of \cite{mestre:2008b} only applies when the number of
distinct population eigenvalues remains finite as matrix dimension
goes to infinity. In practice, this means that the number of distinct
eigenvalues must be negligible with respect to the total number of
eigenvalues. As a further restriction, the number of distinct
eigenvalues and their multiplicities must be known. The only unknown
quantities to be estimated are the locations of the eigenvalues; of
course, this is still a  difficult task. 
\cite{yao:kammoun:najim:2012} propose a more general estimation
procedure that does not require knowledge of the multiplicities, though it
still requires knowledge of the number of distinct population eigenvalues. This setting is
too restrictive for many applications.

The method developed by \cite{karoui:2008} allows for an arbitrary set of population eigenvalues, but does not appear to have good
finite-sample properties.
In fact, our simulations seem to indicate that this estimator
is not even consistent; see Section~\ref{subsub:continuous}.

The first contribution of the present paper is, therefore, to develop an estimator of
the population eigenvalues that is consistent
under large-dimensional asymptotics {\em regardless} of whether or not
they are clustered, and that also performs well in finite sample. This
is achieved through a more precise characterization of the asymptotic
behavior of sample eigenvalues. Whereas existing results only specify
how the eigenvalues behave on average, namely, how many fall in any
given interval, we determine {\em individual} limits.

Our second contribution is to show how this consistent
estimator of population eigenvalues can be used for improved estimation of the covariance matrix when the dimension is large compared to the sample
size. This was already considered in \cite{ledoit:wolf:2012}, but only
in the limited setup where the dimension is smaller than the sample
size. Thanks to the advances introduced in the present paper, we can
also handle the more difficult case where the dimension exceeds the sample
size and the sample covariance matrix is singular.

Our third and final contribution is to show how the same nonlinear
shrinkage formula can be used to estimate the fraction of variation explained by
a given collection of principal components in PCA, which is key in
deciding how many principal components to retain.

The remainder of the paper is organized as follows. 
Section~\ref{sec:asy} presents our estimator of the eigenvalues of the
population covariance matrix under large-dimensional asymptotics.
Section~\ref{sec:sigma} discusses covariance matrix estimation, and
Section~\ref{sec:pca} principal component analysis.
Section~\ref{sec:monte-carlo} studies finite-sample performance via Monte Carlo simulations.
Section~\ref{sec:emp-app} provides a brief empirical application of
PCA to stock return data. Section~\ref{sec:conclusion} concludes.
The proofs of all mathematical results are collected in the appendix.

\section{Estimation of Population Covariance Matrix Eigenvalues}
\label{sec:asy}

\subsection{Large-Dimensional Asymptotics and Basic Framework}
\label{sub:framework}

Let $n$ denote the sample size and $p \defeq p(n)$ the number of
variables. It is assumed that the ratio
$p/n$ converges as $n\to\infty$ to a limit $c\in(0,1) \cup (1,\infty)$ called the {\it concentration}.
The case $c = 1$ is ruled out for
technical reasons. We make the following assumptions.

\begin{itemize} 

\item[(A1)] The population covariance matrix $\Sigma_n$ is a nonrandom
\mbox{$p$-dimensional} positive definite matrix. 

\item[(A2)] $X_n$ is an $n\times p$ matrix of
real independent and identically distributed (i.i.d.)~random
variables with zero mean, unit variance, and finite fourth moment. One only observes $Y_n \defeq
X_n \Sigma_n^{1/2}$, so neither~$X_n$ nor $\Sigma_n$ are
observed on their own.

\item[(A3)] 
$\boldsymbol{\tau}_n\defeq(\tau_{n,1},\ldots,\tau_{n,p})'$ denotes a system of
eigenvalues of~$\Sigma_n$, sorted in increasing order, and 
$(v_{n,1},\ldots,v_{n,p})$ denotes an associated system of eigenvectors. The
empirical distribution function (e.d.f.)\ of the population
eigenvalues is defined as: $\forall t\in\R, \;
H_n(t)\defeq p^{-1}\sum_{i=1}^p\bone_{[\tau_{n,i},+\infty)}(t)$, where~$\bone$
denotes the indicator function of a~set. $H_n$ is called the {\it
  spectral distribution (function)}.
It is assumed that $H_n$ converges weakly to a limit law~$H$,
called the {\it limiting spectral distribution \mbox{(function)}}.

\item[(A4)] $\supp(H)$, the support of $H$, is the union of a
finite number of closed intervals, bounded away from zero and
infinity. Furthermore, there exists a compact interval in $(0, 
\infty)$ that contains $\supp(H_n)$ for all $n$ large enough.

\end{itemize} 

Let $ \boldsymbol{\lambda}_n \defeq (\lambda_{n,1},\ldots,\lambda_{n,p})'$
 denote a
system of eigenvalues of the sample covariance matrix
$S_n \defeq n^{-1} Y_n^\prime Y_n =
n^{-1}\Sigma_n^{1/2}X_n'X_n^{}\Sigma_n^{1/2}$, sorted in increasing order,
and let $(u_{n,1},\ldots,u_{n,p})$ denote an associated system of eigenvectors.
The first subscript, $n$, may be omitted when no confusion is
possible. The e.d.f.\ of the sample eigenvalues is defined as:
$\forall t \in\R, \,
F_n(t)\defeq
p^{-1}\sum_{i=1}^p\bone_{[\lambda_i,+\infty)}(t)$. 
The literature on the eigenvalues of sample covariance matrices under
large-dimensional asymptotics --- also known as {\it random matrix theory} (RMT)
  literature --- is based on a foundational result due to
\cite{marcenko:pastur:1967}. It has been strengthened and broadened by
subsequent authors including \cite{silverstein:1995}, \cite{silverstein:bai:1995},
 \cite{silverstein:choi:1995}, and
\cite{bai:silverstein:1998,bai:silverstein:1999}, among
others. These articles imply that there exists a limiting sample
spectral distribution $F$ such that
\be
\label{eq:convergence}
\forall x\in\R \setminus \{0\}\qquad  F_n(x) \stackrel{\rm
  a.s.}{\longrightarrow} F(x)~.
\ee
In other words, the {\em average} number of sample eigenvalues falling in any given interval is known asymptotically. 

In addition, the existing literature has unearthed important
information about the limiting distribution $F$. \mbox{Silverstein} and Choi (\citeyear{silverstein:choi:1995}) show that $F$ is everywhere continuous except (potentially) at zero, and that the mass that $F$ places at zero is given by
\be
\label{eq:mass0}
F(0)=\max \Bigl \{1-\frac{1}{c}, H(0) \Bigr \}~.
\ee
Furthermore, there is a seminal equation relating $F$ to $H$ and
$c$. Some additional notation is required to present this equation.

For any nondecreasing function $G$ on the real line, 
$m_G$ denotes the {\it Stieltjes transform} of~$G$:
$$
\forall z\in\mathbb{C}^+\qquad m_G(z) \defeq
\int \frac{1}{\lambda-z} \, dG(\lambda)~,
$$
where $\C^+$ denotes the half-plane of complex numbers with strictly
positive imaginary part.

The Stieltjes transform admits a well-known inversion formula:
\be
\label{eq:inversion}
G(b)-G(a)=\lim_{\eta\to0^+}\frac{1}{\pi}\int_a^b\im\bigl
[m_G(\xi+i\eta)\bigr ]d\xi~,
\ee
if $G$ is continuous at $a$ and $b$. Here, and in the remainder of the paper, we shall use the notations $\re(z)$ and
$\im(z)$ for the real and imaginary parts, respectively, of a complex number~$z$, so
that
$$
\forall z \in \mathbb{C} \qquad z = \re(z) + i \cdot \im(z)~.
$$

The most elegant version of the equation relating $F$ to $H$ and
$c$, due to \cite{silverstein:1995}, states that
$m\defeq m_F(z)$ is the unique solution in the set
\be
\label{eq:set}
\left\{m\in\C:-\frac{1-c}{z}+cm\in\C^+\right\}
\ee
to the equation
\be
\label{eq:MP}
\forall z\in\mathbb{C}^+\qquad
m_F(z)=\int \frac{1}{\tau\bigl
  [1-c-c\,z\,m_F(z)\bigr ]-z}\,dH(\tau)~.
\ee

As explained, the Stieltjes transform of $F$, $m_F$, is a function
whose domain is the upper half of the complex plane. It can be
extended to the real line, since
\cite{silverstein:choi:1995} show that:
$\forall\lambda\in\R \setminus \{0\}, \; \lim_{z\in\mathbb{C}^+\to\lambda}m_F(z)\eqdef
\breve{m}_F(\lambda)$ exists. When $c < 1$, $\breve m_F(0)$ also
exists and $F$~has a continuous derivative
$F'=\pi^{-1}\im\left[\breve{m}_F\right]$ on all of $\R$ with $F^\prime
\equiv 0$ on $(-\infty, 0]$. (One should remember that although
the argument of $\breve m_F$ is real-valued now, the output of the
function is still a complex number.) 

For purposes that will become apparent
later, it is useful to reformulate equation \eqref{eq:MP}.
The limiting e.d.f.~of the eigenvalues of $n^{-1} Y_n'Y_n =
n^{-1} \Sigma_n^{1/2} X_n^\prime X_n \Sigma_n^{1/2}$ was defined as
$F$. In addition, define the 
limiting e.d.f.~of the eigenvalues of 
$n^{-1}Y_n Y_n' = n^{-1} X_n\Sigma_nX_n'$ as
$\underline{F}$; note that the eigenvalues of $n^{-1} Y_n'Y_n$ and
$n^{-1} Y_nY_n^\prime$ only differ by $|n-p|$ zero eigenvalues. 
It then holds:
\begin{align}
\forall x\in\R\qquad\underline{F}(x)&=(1-c)\,
\bone_{[0,\infty)}(x)+c\,F(x) \label{eq:mF1} \\
\forall x\in\R\qquad
F(x)&=\frac{c-1}{c}\bone_{[0,\infty)}(x)+\frac{1}{c}\,\underline{F}(x)
\label{eq:mF2} \\
\forall z \in \C^+ \qquad m_{\underline{F}}(z)&=\frac{c-1}{z}+c\,
m_F(z) \label{eq:mF3}\\
\forall z \in \C^+ \qquad m_F(z)&=\frac{1-c}{c\, z}+\frac{1}{c} \, m_{\underline{F}}(z)~.
\label{eq:mF4}
\end{align}
(Recall here that $F$ has mass $(c-1)/c$ at zero when $c > 1$,
so that both $F$ and $\underline F$ are nonnegative functions
indeed for any value $c > 0$.)

With this notation, equation (1.13) of \cite{marcenko:pastur:1967}
reframes equation~\eqref{eq:MP} as: for each $z\in\C^+$, $m\defeq m_{\underline{F}}(z)$ is the unique solution in $\C^+$ to the equation
\be
\label{eq:MP-bar}
m=- \left [
{z-c\int \frac{\tau}{1+\tau \,
    m} dH(\tau)} \right ]^{-1}.
\ee
% It is clear from these equations that consistent estimation of $F$ and
% $m_{F}$, respectively, immediately results in consistent estimation of
% of $\underline F$ and $m_{\underline{F}}$, respectively, as well; of
% course, in doing so, one uses the consistent-by-definition estimator $\widehat c_n \defeq
% p/n$ for the concentration $c$.
While in the case $c < 1$, $\breve m_F(0)$ exists and $F$
is continuously differentiable on all of $\R$, as mentioned above, 
in the case $c > 1$, $\breve m_{\underline F}(0)$ exists and $\underline
F$ is continuously differentiable on all of~$\R$.

\subsection{Individual Behavior of Sample Eigenvalues: the QuEST Function}
\label{sub:quest}

We introduce a nonrandom multivariate function called the {\it Quantized Eigenvalues Sampling Transform},
or QuEST for short, which discretizes, or {\em
  quantizes},  the relationship between $F$, $H$, and $c$ defined in equations \eqref{eq:convergence}--\eqref{eq:inversion}. For any positive integers $n$ and $p$,
the QuEST function, denoted by $Q_{n,p}$, is defined as 
\begin{align}
Q_{n,p}:[0,\infty)^p&\longrightarrow[0,\infty)^p \label{eq:quest1}\\
\mathbf{t}\defeq\left(t_1,\dots,t_p\right)^\prime&\longmapsto
Q_{n,p}(\mathbf{t})\defeq\left(q_{n,p}^1(\mathbf{t}),\ldots,q_{n,p}^p(\mathbf{t})\right)^\prime~,\label{eq:quest2} 
\end{align}
where 
$\forall z\in\C^+\quad m\defeq m_{n,p}^\mathbf{t}(z)$ is the unique solution in the set 
\be
\left\{m\in\C:-\frac{n-p}{nz}+\frac{p}{n}\,m\in\C^+\right\}\label{eq:questset}
\ee
to the equation 
\be
m=\frac{1}{p}\sum_{i=1}^p\frac{1}{\displaystyle t_i\left(1-\frac{p}{n}-\frac{p}{n}\,z\,m\right)-z}~,\label{eq:questMP}
\ee
\begin{align}
\forall x\in\R\qquad F_{n,p}^\mathbf{t}(x) &\defeq 
\begin{cases}
\displaystyle\max\left\{1-\frac{n}{p},\frac{1}{p}\sum_{i=1}^p\bone_{\{t_i=0\}} \right\}&\text{if $x=0$~,}\\
\displaystyle\lim_{\eta\to0^+}\frac{1}{\pi}\int_{-\infty}^x
\im\left[m_{n,p}^\mathbf{t}(\xi+i\eta)\right]\,d\xi & \text{otherwise~,}
\end{cases}\label{eq:questF}\\
\forall u\in[0,1]\qquad \left(F_{n,p}^\mathbf{t}\right)^{-1}(u)&\defeq \sup\{x\in\R:F_{n,p}^\mathbf{t}(x)\leq u\}~,\label{eq:inverse}\\
\mbox{and}\qquad\forall i=1,\ldots,p\qquad 
q_{n,p}^i(\mathbf{t})&\defeq
p\displaystyle\int_{(i-1)/p}^{i/p}\left(F_{n,p}^\mathbf{t}\right)^{-1}(u)\,du~.\label{eq:condexp}
\end{align}
It is obvious 
that equation (\ref{eq:questset}) quantizes equation (\ref{eq:set}), 
that equation (\ref{eq:questMP}) quantizes equation (\ref{eq:MP}),
and that equation~(\ref{eq:questF}) quantizes
equations (\ref{eq:mass0}) and (\ref{eq:inversion}). 
Thus,
$F_{n,p}^\mathbf{t}$ is the limiting distribution (function) of sample
eigenvalues corresponding to the population spectral distribution
(function) $p^{-1}\sum_{i=1}^p\bone_{[t_{i},+\infty)}$. Furthermore,
by equation (\ref{eq:inverse}), $\left(F_{n,p}^\mathbf{t}\right)^{-1}$
represents the inverse spectral distribution function, also known as
the {\em quantile} function.

\begin{remark}[Definition of Quantiles]\label{rem:quant}
\rm
The standard definition of the $(i-0.5)/p$ quantile of~$F_{n,p}^\mathbf{t}$ is
$\left(F_{n,p}^\mathbf{t}\right)^{-1}((i-0.5)/p)$, where
$\left(F_{n,p}^\mathbf{t}\right)^{-1}$ is defined in equation
(\ref{eq:inverse}). It turns out, however, that the `smoothed' version
$q_{n,p}^i(\mathbf{t})$ given in equation (\ref{eq:condexp}) leads to
improved accuracy, higher stability, and faster
computations of our numerical algorithm, to be detailed below, in practice.

Since $F_n$ is an empirical distribution (function), its quantiles are not
uniquely defined. For example, the statistical software {\sc R}
offers nine different versions of sample quantiles in its function {\tt
  quantile}; version 5 corresponds to our convention of considering
$\lambda_{n,i}$ as the $(i-0.5)/p$ quantile of $F_n$.~\qed
\end{remark}
Consequently, a set of $(i-0.5)/p$ quantiles ($i=1, \ldots, p$) is
given by $Q_{n,p}(\tf)$ for $F_{n,p}^\mathbf{t}$ and  is given by
$\boldsymbol{\lambda}_n$ for $F_n$. 
The relationship between
$Q_{n,p}(\tf)$ and $\boldsymbol{\lambda}_n$ is further elucidated by
the following theorem.
\begin{theorem}
\label{theo:individual}
If Assumptions (A1)--(A4) are satisfied, then
\be
\frac{1}{p}\sum_{i=1}^p\left[q_{n,p}^i(\boldsymbol{\tau}_n)-\lambda_{n,i}\right]^2\stackrel{\rm
  a.s.}{\longrightarrow}0~.\label{eq:individual}
\ee
\end{theorem}
Theorem \ref{theo:individual} states that the sample eigenvalues
converge individually to their nonrandom QuEST function
counterparts. This individual notion of convergence is defined as the
Euclidian distance between the vectors $\boldsymbol{\lambda}_n$ and
$Q_{n,p}(\boldsymbol{\tau}_n)$, normalized by the matrix dimension~$p$.
It is the appropriate normalization because, as $p$ goes to
infinity, the left-hand side of equation~(\ref{eq:individual})
approximates the $L^2$ distance between the functions $F_n^{-1}$ and
$\left(F_{n,p}^{\boldsymbol{\tau}_n}\right)^{-1}$. This metric can be
thought of as a `cross-sectional' mean squared error, in the same way
that $F_n$ is a cross-sectional distribution function.

Theorem \ref{theo:individual} improves over the well-known results from the random matrix theory literature reviewed in Section \ref{sub:framework} in two significant ways.
\begin{enumerate}
\renewcommand{\labelenumi}{\arabic{enumi})}
\item It is based on the $p$ population eigenvalues
  $\boldsymbol{\tau}_n$, not the limiting spectral distribution
  $H$. Dealing with $\boldsymbol{\tau}_n$ (or, equivalently, $H_n$) is
  straightforward because it is integral to the actual data-generating
  process; whereas dealing with $H$ is more delicate because we do not
  know how $H_n$ converges to $H$. Also there are potentially
  different $H$'s that $H_n$ could converge to, depending on what we
  assume will happen as the dimension increases. 

\item Theorem \ref{theo:individual} characterizes the {\em individual}
  behavior of the sample eigenvalues, whereas equation
  (\ref{eq:convergence}) only characterizes their {\em average}
  behavior, namely, what proportion falls in any given
  interval. Individual results are more precise than average
  results. Thus, Theorem \ref{theo:individual} shows that the sample
  eigenvalues are better behaved under large-dimensional asymptotics
  than previously thought. 

\end{enumerate}
Both of these improvements are made possible thanks to the introduction of the QuEST function.
In spite of the apparent complexity of the mathematical definition of the QuEST function, 
it can be computed quickly and efficiently along with 
its analytical Jacobian as evidenced by Figure \ref{fig:speed}, and it behaves well numerically.
\begin{center}
\captionsetup{type=figure}
\includegraphics[scale=\scale]{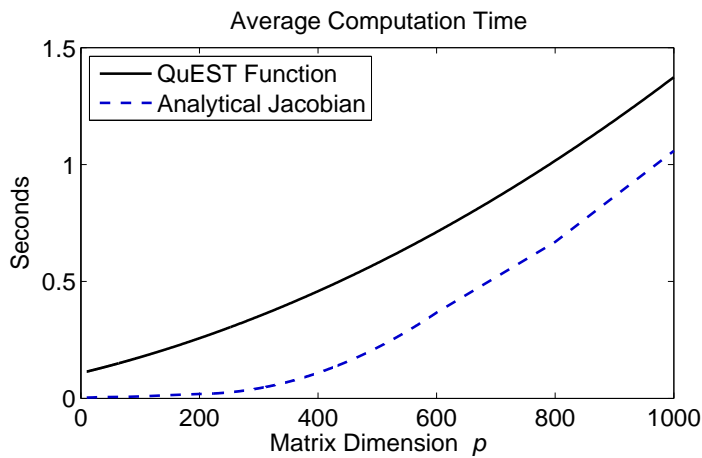}
\captionof{figure}{Average computation time for the QuEST function and its analytical Jacobian. 
The setup is the same as in Figure \ref{fig:conv1} below. 
The QuEST function and its analytical Jacobian are programmed in Matlab. 
The computer is a $2.4$ GHz desktop Mac.}\label{fig:speed}
\end{center}
% software references for Olivier [do not delete this comment]:
% comptime01.eps, optim19, QuESTimer_plotter02.m
% tau51_*.mat, QuESTimer03.m, 15 Jul 2013

\subsection{Consistent Estimation of Population Eigenvalues}
\label{sub:estimator}

Once the truth of Theorem \ref{theo:individual} has been established,
it becomes tempting to construct an estimator of population covariance
matrix eigenvalues simply by minimizing the expression on the
left-hand side of equation (\ref{eq:individual}) over all possible
sets of population eigenvectors. This is exactly what we do in the
following theorem.
\begin{theorem}
\label{theo:estimator}
Suppose that Assumptions (A1)--(A4) are satisfied. Define 
\begin{equation} \label{eq:optim}
\widehat{\boldsymbol{\tau}}_n \defeq \argmin_{\mathbf{t}\in[0,\infty)^p}
\,\frac{1}{p}\sum_{i=1}^p\left[q_{n,p}^i(\mathbf{t})-\lambda_{n,i}\right]^2,
\end{equation} 
where $\boldsymbol{\lambda}_n \defeq (\lambda_{n,1}, \ldots, \lambda_{n,p})^\prime$
are the sample covariance matrix eigenvalues, and $Q_{n,p}(\mathbf{t})\defeq \left(q_{n,p}^1(\mathbf{t}), \ldots, q_{n,p}^p(\mathbf{t})\right)^\prime$ is the nonrandom QuEST function
defined in equations \eqref{eq:quest1}--\eqref{eq:questMP}; both
$\widehat{\boldsymbol{\tau}}_n$ and $\boldsymbol{\lambda}_n$  are
assumed sorted in increasing order. Let $\widehat{\tau}_{n,i}$ denote
the $i$th entry of $\widehat{\boldsymbol{\tau}}_n$ $(i=1,\ldots,p)$, and let
$\boldsymbol{\tau}_n\defeq(\tau_{n,1},\ldots,\tau_{n,p})'$ denote the
population covariance matrix eigenvalues sorted in increasing
order. Then
\be
\frac{1}{p}\sum_{i=1}^p\left[\widehat{\tau}_{n,i}-\tau_{n,i}\right]^2\stackrel{\rm
  a.s.}{\longrightarrow}0~.\label{eq:consistent}
\ee
\end{theorem}
Theorem \ref{theo:estimator} shows that the estimated eigenvalues
converge {\em individually} to the population eigenvalues, in the same
sense as above, using the dimension-normalized Euclidian distance.

\begin{remark}
\rm
Mathematically speaking, equation (\ref{eq:optim}) performs two tasks:
it projects $\boldsymbol{\lambda}_n$ onto the image of the QuEST
function, and then inverts the QuEST function. Since the image of the
QuEST function is a strict subset of $[0,\infty)^p$,
$\boldsymbol{\lambda}_n$ will generally be outside of it. It is the
first of these two tasks that gets around any potential ill-posedness
by {\em regularizing} the set of observed sample eigenvalues. 

Practically speaking, both tasks are performed simultaneously by a nonlinear optimizer. We use a standard off-the-shelf commercial software called SNOPT\textsuperscript{\texttrademark}
 (\mbox{Version 7.4}), see \cite{gill:murray:saunders:2002},
but other choices may work well too.~\qed
\end{remark}

\subsection{Comparison with Other Approaches}

\cite{karoui:2008} also attempts to discretize equation (\ref{eq:MP}) and invert it, but he opts for a completely
opposite method of discretization which does not exploit the natural discreteness of the population spectral distribution for finite $p$. If the population
spectral distribution $H_n$ is approximated by a convex linear combination of step functions
\begin{equation*}\label{eq:H-tilde}
\forall x\in\R\quad \widetilde{H}(x) \defeq \sum_{i=1}^pw_i\bone_{\{x\geq t_i\}}
\quad\mbox{where}\quad\forall i=1,\ldots,p\quad
t_i\geq0,\;w_i\geq0\;\mbox{, and}\;\sum_{i=1}^pw_i=1~,
\end{equation*} 
then in the optimization problem (\ref{eq:optim}), we keep the weights
$w_i$ ($i=1,\ldots,p$) fixed at $1/p$ while varying the location
parameters $t_i$ ($i=1,\ldots,p$). In contrast, \cite{karoui:2008}
does exactly the reverse:
he keeps the location parameters~$t_i$ fixed on a grid while varying
the weights~$w_i$. Thus, \cite{karoui:2008} projects the population
spectral distribution onto a ``dictionary''. Furthermore, instead of matching population
eigenvalues to sample eigenvalues on $\R$, he matches a
function of $m_{\widetilde H}$ to a function of $m_{F_n}$ on $\C^+$, which makes his algorithm
relatively complicated; see \citeauthor{ledoit:wolf:2012} (2012, pages
1043--1044).
Despite our best efforts, we were unable to
replicate his convergence results in Monte Carlo simulations: in our
implementation, his estimator performs poorly overall and does not
even appear to be consistent; see Section~\ref{subsub:continuous}.
Unless someone circulates an implementation
of the algorithm described in \cite{karoui:2008} that works, we have to write
off this approach as impractical.

Another related paper is the one by \cite{ledoit:wolf:2012}. They use
the same discretization strategy as \cite{karoui:2008} (fix location
parameters and vary weights) but, as we do here, match population
eigenvalues to sample eigenvalues on the real line. Unlike we do here,
they measure closeness by a sup-distance
rather than by the Euclidean distance.
 \cite{ledoit:wolf:2012} only consider the case $p<n$. 
Unfortunately, their nonlinear
optimizer no longer converges reliably in the case $p > n$, as we found
out in subsequent experiments; this necessitated the development of
the alternative discretization strategy described above, as well as
the change from sup-distance to Euclidean distance to measure closeness.

Furthermore,
\cite{ledoit:wolf:2012} are not directly interested in estimating 
the population eigenvalues;
it is just an intermediary step towards their ultimate
objective, which is the estimation of the covariance matrix itself.
Therefore they do not report any Monte Carlo simulations on the finite-sample
behavior of their estimator of the population eigenvalues. 

In any case, the aim
of the present paper is to develop an estimator of the population eigenvalues that works also when $p > n$,
so the approach of \cite{ledoit:wolf:2012} is ruled out. 
The different discretization
strategy that we employ here, together with the alternative distance
measure, enables us to construct an estimator of
$\boldsymbol{\tau}_n$ that works across both cases $p<n$ and $p>n$.
It is important to point out that the new estimator of population eigenvalues is not only
more general, in the sense that it also works for the case $p>n$, but
it also works better for the case $p<n$; see
Section~\ref{ss:mc-sigma}.

As for the papers of \cite{mestre:2008b} and
\cite{yao:kammoun:najim:2012}, their methods are based on contour
integration of analytic functions in the complex plane. They can only
extract a finite number~$\bar{M}$ of functionals of $H_n$, such as the
locations of high-multiplicity eigenvalue clusters, or the trace of
powers of $\Sigma_n$. The main difference with our method is that we
extract many more items of information: namely, $p$ population
eigenvalues. This distinction is crucial because the \mbox{ratio $\bar{M}/p$}
vanishes asymptotically. It explains why we are able to recover the
whole population spectrum in the general case, whereas they are not.

\section{Covariance Matrix Estimation}
\label{sec:sigma}

The estimation of the covariance matrix $\Sigma_n$ is already considered by
\cite{ledoit:wolf:2012}, but only for the case $p < n$. In particular,
they propose a nonlinear shrinkage approach, which we will now extend
to the case $p > n$. To save space, the reader is referred to their
paper for a more detailed discussion of the nonlinear shrinkage
methodology and a comparison to other estimation strategies of
large-dimensional covariance matrices, such as the linear shrinkage
estimator of \cite{ledoit:wolf:2004a}.

\subsection{Oracle Shrinkage}
\label{ss:oracle}

The starting point is to restrict attention to {\it
  rotation-equivariant} estimators of $\Sigma_n$. 
 To~be more specific, let $W$ be an arbitrary
$p$-dimensional rotation matrix. Let $\widehat \Sigma_n \defeq \widehat \Sigma_n(Y_n)$
be an estimator of $\Sigma_n$. Then the
estimator is said to be {\it rotation-equivariant} if 
it satisfies $\widehat \Sigma_n(Y_n W) = W^\prime \widehat \Sigma_n(Y_n) W$. 
In other words, the estimate based on the rotated data equals the
rotation of the estimate based on the original data.
In the absence of any {\it a priori} knowledge about the structure of
$\Sigma_n$, such as sparseness or a factor model, it is natural to
consider only estimators of $\Sigma_n$ that are rotation-equivariant.

The class of rotation-equivariant estimators of the covariance that
are a function of the sample covariance matrix
is constituted of all the estimators that have the same eigenvectors
as the sample covariance matrix; for example, see
\citeauthor{perlman:2007} (2007, Section~5.4). 
Every such rotation-equivariant estimator is thus of
the form
\begin{equation} 
\label{eq:rot}
U_n D_n U_n^\prime \quad \mbox{where} \quad
D_n \defeq \diag(d_1, \ldots, d_p) \mbox{ is diagonal}~,
\end{equation} 
and where $U_n$ is the matrix whose $i$th column is the sample
eigenvector $u_i \defeq u_{n,i}$. This is the class of
rotation-equivariant estimators already studied by
\cite{stein:1975,stein:1986}.

We can rewrite the expression for such a
rotation-equivariant estimator as
\begin{equation} \label{eq:rot-equ}
U_n D_n U_n' = \sum_{i=1}^p d_i \cdot u_i^{} u_i^\prime~.
\end{equation} 
This alternative expression shows that any such rotation-equivariant
estimator is a linear combination of $p$ rank-1 matrices $u_i^{}
u_i^\prime \; (i=1, \ldots, p)$. But since the $\{u_i\}$ form an orthonormal basis in~$\R^p$,
the resulting estimator is still of full rank $p$, provided that all the
weights $d_i \; (i = 1, \ldots, p)$ are strictly positive.

\begin{remark}[Rotation-equivariant Estimators versus Structured Estimators]
\label{rem:rot}
\rm
By construction, the class \eqref{eq:rot} of rotation-invariant
estimators have the same eigenvectors as the sample covariance
matrix. In particular, consistent estimation of the
covariance matrix is not possible under large-dimensional asymptotics.

Another approach would be to impose additional structure on the
estimation problem, such as sparseness \citep{bickel:levina:2008}, a
graph model \citep{bala:et:al:2008}, or an (approximate) factor model 
\citep{fan:fan:lv:2008}.\footnote{We only give one representative
  reference for each field here to save space.} 
The advantage of doing so is that, under
suitable regularity conditions, consistent estimation of the
covariance matrix is possible. The disadvantage is that if the assumed
structure is misspecified, the estimator of the covariance matrix can
be arbitrarily bad; and whether the structure is correctly specified
can never be verified from the data alone.

Rotation-equivariant estimators are widely and successfully used in
practice in situations where knowledge on additional structure is not
available (or doubtful).
This is evidenced by the many citations to
\cite{ledoit:wolf:2004a} who propose a linear shrinkage estimator that
also belongs to the class~\eqref{eq:rot}; for example, see the
beginning of Section~\ref{ss:mc-sigma}. Therefore, developing a new, nonlinear
shrinkage estimator that outperforms this previous proposal will
be of substantial interest to applied researchers in our opinion.~\qed
\end{remark}

The first objective is to find the matrix in the class \eqref{eq:rot} of
rotation-equivariant estimators that is
closest to $\Sigma_n$. To measure distance, we choose the
Frobenius norm defined as
\begin{equation} \label{e:fro}
||A||_F \defeq \sqrt{\tr(A
  A^\prime)/r} \quad \mbox{ for any matrix $A$ of dimension $r \times m$}~.
\end{equation} 
(Dividing by the dimension of the square matrix $A A^\prime$ inside
the root is not
standard, but we do this for asymptotic purposes so that the Frobenius
norm remains constant equal to one for the identity matrix regardless of the
dimension; see \cite{ledoit:wolf:2004a}.)
As a result, we end up with the following
minimization problem:
$$
\min_{D_n} ||U_n D_n U_n^\prime - \Sigma_n||_F~.
$$
Elementary matrix algebra shows that its solution is
\begin{equation}\label{e:star}
D_n^* \defeq \diag(d_1^*, \ldots, d_p^*) \quad
\mbox{where} \quad 
d_i^* \defeq u_i^\prime \Sigma_n u_i^{}\; \mbox{ for } i = 1, \ldots, p~.
\end{equation}

Let $y \in \R^p$ be a random vector with covariance matrix $\Sigma_n$, drawn independently from the sample covariance matrix $S_n$. We can think of $y$ as an out-of-sample observation.
Then $d_i^*$ is recognized as the variance of
the linear combination $u_i^\prime y$, conditional on $S_n$. In view of the expression~\eqref{eq:rot-equ}, it makes
intuitive sense that the matrices $u_i^{} u_i^\prime$ whose associated
linear combination $u_i^\prime y$ have higher out-of-sample variance should receive higher weight in
computing the estimator of $\Sigma_n$.

The finite-sample optimal estimator is thus given by
\begin{equation}\label{e:optimal-fs}
S_n^* \defeq U_n D_n^* U_n^\prime \quad \mbox{where} \quad
D_n^* \mbox{ is defined as in \eqref{e:star}}~.
\end{equation}

Clearly $S_n^*$ is not a {\em feasible} estimator because it depends on knowing the population covariance matrix.
By generalizing the Mar\v{c}enko-Pastur equation (\ref{eq:MP}),
\cite{ledoit:peche:2011} show that~$d_i^*$ can be approximated by the
asymptotic quantities
\begin{equation}\label{eq:oracle0}
d_i^{or} \defeq 
\left \{
\renewcommand{\arraystretch}{1.5}
\begin{array}{cc}
\displaystyle\frac{1}{(c-1) \, \breve m_{\underbar{\mbox{\tiny $F$}}} (0)} ~, &
    \mbox{ if } \lambda_i = 0 \mbox{ and } c>1\\
\displaystyle\frac{\lambda_i}{\bigl |1 - c -c \, \lambda_i \, \breve
  m_F(\lambda_i) \bigr |^2}~,& \mbox{ otherwise } 
\end{array}
\right .
\quad \mbox{ for } i = 1, \ldots, p~,
\end{equation}
from which they deduce their oracle estimator
\begin{equation}\label{eq:oracle}
S_n^{or} \defeq U_n D_n^{or} U_n^\prime \quad 
\mbox{where} \quad 
D_n^{or} \defeq \diag(d_1^{or}, \ldots, d_p^{or})~.
\end{equation} 

The key difference between $D^*_n$ and $D^{or}_n$ is that the former
depends on the unobservable population covariance matrix, whereas the
latter depends on the limiting distribution of sample eigenvalues, $F$,
which makes it amenable to consistent estimation. It turns out that this
estimation problem is solved if a consistent estimator of the
population eigenvalues $\boldsymbol{\tau}_n$ is available.

\subsection{Nonlinear Shrinkage Estimator}
\label{ss:nonlin}

\subsubsection{The Case \texorpdfstring{$p < n$}{p < n}}

We start with the case $p < n$, which was already considered by \cite{ledoit:wolf:2012}.
\cite{silverstein:choi:1995} show how
the support of $F$, denoted by $\supp(F)$, is determined; also see
Section~2.3 of \cite{ledoit:wolf:2012}. $\supp(F)$ is
seen to be the union of a finite number of disjoint compact intervals,
bounded away from zero. To simplify the discussion, we will assume
from here on that $\supp(F)$ is a single compact interval, bounded
away from zero, with $F^\prime > 0$ in the interior of this interval. 
But if $\supp(F)$ is the union of a finite
number of such intervals, the arguments presented in this section as
well as in the remainder of the paper apply separately to each
interval. In particular, our consistency results presented below
can be easily extended to this more general case.

Recall that, for any $\mathbf{t}\defeq(t_1,\ldots,t_p)'\in[0,+\infty)^p$, equations (\ref{eq:questset})--(\ref{eq:questMP}) define $m_{n,p}^\mathbf{t}$ as the Stieltjes transform of $F_{n,p}^\mathbf{t}$, the limiting distribution of sample eigenvalues corresponding to the population spectral distribution $p^{-1}\sum_{i=1}^p\bone_{[t_{i},+\infty)}$. Its domain is the strict upper half of the complex plane, but it can be
extended to the real line since
\cite{silverstein:choi:1995} prove that
$\forall\lambda\in\R-\{0\}\quad
\lim_{z\in\mathbb{C}^+\to\lambda}m_{n,p}^\mathbf{t}(z) \eqdef
\breve{m}_{n,p}^\mathbf{t}(\lambda)$ exists. 

\cite{ledoit:wolf:2012} show how a consistent estimator of $\breve
m_F$ can be derived from a consistent estimator of $\boldsymbol{\tau}_n$, such as $\widehat{\boldsymbol{\tau}}_n$ defined in Theorem \ref{theo:estimator}. Their Proposition~4.3 
establishes that $\breve{m}_{n,p}^{\widehat{\boldsymbol{\tau}}_n}(\lambda) \to \breve
m_F(\lambda)$ uniformly in $\lambda \in \supp(F)$, except for two
arbitrarily small regions at the lower and upper end of
$\supp(F)$. Replacing $\breve m_F$ with $\breve{m}_{n,p}^{\widehat{\boldsymbol{\tau}}_n}$
and $c$ with $p/n$ in \cite{ledoit:peche:2011}'s oracle
quantities $d_i^{or}$ of \eqref{eq:oracle0} yields 
\begin{equation}\label{eq:nonlin0-p-less-n}
\wh d_i \defeq \frac{\lambda_i}
 {\left|1 - \displaystyle\frac{p}{n} -  \frac{p}{n} \, \lambda_i \cdot 
  \breve{m}_{n,p}^{\widehat{\boldsymbol{\tau}}_n}(\lambda_i) \right|^2}
\quad \mbox{ for } i = 1, \ldots, p~.
\end{equation}
(Note here that in the case $p < n$, all sample eigenvalues $\lambda_i$ are positive
almost surely, for $n$ large enough, by the results of \cite{bai:silverstein:1998}.)
In turn,  the {\it bona fide}
nonlinear shrinkage estimator of $\Sigma_n$ is obtained as:
\begin{equation}\label{eq:nonlin}
\wh S_n \defeq U_n \wh D_n U_n^\prime \quad 
\mbox{where} \quad 
\wh D_n \defeq \diag(\wh d_1, \ldots, \wh d_p)~.
\end{equation} 

\subsubsection{The Case \texorpdfstring{$p > n$}{p > n}}

We move on to the case $p > n$, which was not considered by
\cite{ledoit:wolf:2012}. In this case, $F$ is a mixture distribution
with a discrete part and a continuous part. 
The discrete part is a
point mass at zero with mass $(c-1)/c$. The continuous part has total mass
$1/c$ and its support is the union of a finite number of disjoint
intervals, bounded away from zero; again, see
\cite{silverstein:choi:1995}.

It can be seen from equations \eqref{eq:mF1}--\eqref{eq:mF4}
that $\underline F$ corresponds to the continuous part of~$F$, scaled
to be a proper distribution (function): $\lim_{t \to \infty}
\underline F(t) = 1$. Consequently, $\supp(F) = \{0\}
\cup \supp(\underline F)$. To simplify the discussion, we will assume
from here on that $\supp(\underline F)$ is a single compact interval, bounded
away from zero, with $\underline F^\prime > 0$ in the interior of this interval. 
But if $\supp(\underline F)$ is the union of a finite
number of such intervals, the arguments presented in this section as
well as in the remainder of the paper apply separately to each
interval. In particular, our consistency results presented below
can be easily extended to this more general case.

The oracle quantities $d_i^{or}$ of \eqref{eq:oracle0} involve
$\breve m_{\underline F} (0)$ and $\breve m_F(\lambda_i)$ for various
$\lambda_i > 0$; recall that~$\breve m_{\underline F} (0)$ exists in
the case $c > 1$.
% From the definition of the Stieltjes transform
% \eqref{eq:stieltjes} it can be seen that 
% \begin{equation}\label{eq:1/X}
% \breve m_{\underline F} (0) = \E(X^{-1}) \mbox{ with } X \sim \underline{F}~.
% \end{equation} 
% Since the support of $\underline F$ is bounded away from zero and
% infinity, $\breve m_{\underline F} (0) \in (0, \infty)$ is well-defined.
% Let $j$ denote the smallest integer for which $\lambda_j > 0$. Then a
% consistent estimator of $\breve m_{\underline F} (0)$ is given by
% \begin{equation} \label{eq:m-0}
% \widehat{\breve m_{\underline F} (0)} \defeq
% \frac{\widehat c_n}{n-j+1} \sum_{i=j}^n \frac{1}{\lambda_i}~.
% \end{equation} 
% A consistent estimator of $\breve m_{\underline F} (0)$ is given by the
% unique solution in $(0, + \infty)$ to the equation
% \begin{equation} 
% \frac{1}{p}\sum_{i=1}^p\frac{1}{1+\widehat{t}_i \, \breve m_{\underline{F}}(0)} =
%   1-\frac{n}{p}~,
% \end{equation} 
% where $\hat \tf_n \defeq \widehat \tf \defeq (\wh t_1, \ldots, \wh
% t_p)^\prime$ is defined as in Theorem~\ref{theo:FH-hat}.

% {\bf Dot: the above is what you supplied to me by e-mail. 
% However, I think it should be the one below \ldots }

Using the original \MP equation \eqref{eq:MP-bar},
a strongly consistent estimator of the quantity~$\breve m_{\underline F} (0)$ is the unique solution $m\defeq\wh{\breve m_{\underbar{\mbox{\tiny $F$}}}(0)}$ in $(0, \infty)$ to the equation
\begin{equation} 
m = \left [
\frac{1}{n} \sum_{i=1}^p \frac{\wh \tau_i}{1 + \wh \tau_i \, m}
\right ]^{-1} \! \! ,
\end{equation} 
where $\widehat{\boldsymbol{\tau}}_n \defeq (\wh \tau_1, \ldots, \wh
\tau_p)^\prime$ is defined as in Theorem~\ref{theo:estimator}.

Again, since $\widehat{\boldsymbol{\tau}}_n$ is consistent for $\boldsymbol{\tau}_n$, Proposition 4.3 of \cite{ledoit:wolf:2012} implies
that $\breve{m}_{n,p}^{\widehat{\boldsymbol{\tau}}_n}(\lambda) \to \breve
m_F(\lambda)$ uniformly in $\lambda \in \supp(\underline F)$, except for two
arbitrarily small regions at the lower and upper end of
$\supp(\underline F)$.

Finally, the {\it bona fide} nonlinear shrinkage estimator of
$\Sigma_n$ is obtained as
\eqref{eq:nonlin} but now with \begin{equation}\label{eq:nonlin0-p-greater-n}
\wh d_i \defeq 
\left \{
\begin{array}{cc}\displaystyle
\frac{\lambda_i}
 {\left|1 - \frac{p}{n} -  \frac{p}{n} \, \lambda_i \cdot 
  \breve{m}_{n,p}^{\widehat{\boldsymbol{\tau}}_n}(\lambda_i) \right|^2}~,& \mbox{ if } \lambda_i > 0 \\
\displaystyle\frac{1}{\left(\frac{p}{n}-1\right) \, \wh{\breve m_{\underbar{\mbox{\tiny $F$}}} (0)}}~, &
    \mbox{ if } \lambda_i = 0
\end{array}
\right .
\quad \mbox{ for } i = 1, \ldots, p~,
\end{equation}

\subsection{Strong Consistency}
\label{ss:nonlin-cons}

The following theorem establishes that our nonlinear shrinkage
estimator, based on the estimator $\widehat{\boldsymbol{\tau}}_n$ of Theorem~\ref{theo:estimator},
is strongly consistent for the oracle estimator across both cases $p<n$ and $p> n$.

\begin{theorem}\label{prop:oracle-cons} 

Let $\widehat{\boldsymbol{\tau}}_n$ be an estimator of the
eigenvalues of the population covariance matrix satisfying
$p^{-1}\sum_{i=1}^p\left[\widehat{\tau}_{n,i}-\tau_{n,i}\right]^2
\toas 0$.
Define the nonlinear shrinkage estimator $\wh S_n$ as
in~\eqref{eq:nonlin}, where the $\wh d_i$ are as in
\eqref{eq:nonlin0-p-less-n} in the case $p< n$ and as in
\eqref{eq:nonlin0-p-greater-n} in the case $p > n$.

\smallskip
Then $||\widehat S_n - S_n^{or}||_F \toas 0$.
\end{theorem}

\begin{remark}\label{rem:c=1}
\rm
We have to rule out the case $c=1$ (or $p = n$) for mathematical
reasons. 

First, we need $\supp(\underline F)$ to be bounded away from
zero to establish various consistency results. 
But when $c=1$, then $\supp(\underline F)$ can start at zero, that is, there exists $u > 0$
such that $F^\prime(\lambda) > 0$ for all $\lambda \in (0, u)$. 
This was already established by \cite{marcenko:pastur:1967}
for the special case when $H$ is a point mass at one. In particular,
the resulting (standard) \MP distribution $F$ has density function
$$
F^\prime(\lambda) = 
\left \{
\begin{array}{ll}
\frac{1}{2\pi \lambda c}\sqrt{(b-\lambda)(\lambda-a)}~, & \mbox{ if $a \le \lambda \le b$}~, \\
0~, & \mbox{ otherwise}~,
\end{array}
\right .
$$
and has point mass $(c-1)/c$ at the origin if $c > 1$, where $a \defeq (1 -
\sqrt{c})^2$ and $b \defeq (1 + \sqrt{y})^2$; for example, see
\citeauthor{bai:silverstein:2010} (2010, Section~3.3.1).

Second, we also need the assumption $c \neq 1$ `directly' in the proof
of Theorem~\ref{prop:oracle-cons} to demonstrate that the summand
$D_1$ in \eqref{eq:d1-d2} converges to zero.

Although the case $c=1$ is not covered by the mathematical treatment,
we can still address it in Monte Carlo simulations; see
Section~\ref{ss:mc-sigma}.~\qed
\end{remark}

\section{Principal Component Analysis}
\label{sec:pca}

Principal component analysis (PCA) is one of the oldest and best-known
techniques of multivariate analysis, dating back to \cite{pearson:1901}
and \cite{hotelling:1933}; for a comprehensive treatment, see
\cite{jolliffe:2002}.

\subsection{The Central Idea and the Common Practice}
\label{ss:pca-1}

The central idea of PCA is to reduce the dimensionality of a data set
consisting of a large number of interrelated variables, while
retaining as much as possible of the variation present in the data
set. This is achieved by transforming the original variables to a new set
of {\it uncorrelated} variables, the principal components, which
are ordered so that the `largest' {\it few} retain most of the variation
present in {\it all} of the original variables.

More specifically, let $y \in \R^p$ be a random vector
with covariance matrix $\Sigma$; in this section, it will be
convenient to drop the subscript $n$ from the covariance
matrix and related quantities. Let
$((\tau_{1},\ldots,\tau_{p});(v_{1},\ldots,v_{p}))$ denote a system of
eigenvalues and eigenvectors of~$\Sigma$. To be consistent with our
former notation, we assume that the eigenvalues $\tau_i$ are sorted in
{\it increasing} order. Then the principal components of $y$ are given
by $v_1^\prime y, \ldots, v_p^\prime y$. Since the eigenvalues
$\tau_i$ are sorted in increasing order, the principal component with
the largest variance is  $v_p^\prime y$ and the principal
component with the smallest variance is 
$v_1^\prime y$. The eigenvector $v_i$ is called the {\it
  vector of coefficients} or {\it loadings} for the $i$th principal
component ($i = 1, \ldots, p$).

Two brief remarks are in order. First, some authors use the term
{\it principal components} for the eigenvectors $v_i$; but we agree with
\citeauthor{jolliffe:2002} (2002, Section~1.1) that this usage is confusing and that it is
preferable to reserve the term for the derived variables $v_i^\prime y$.
Second, in the PCA literature, in contrast to the bulk of the
multivariate statistics literature, the eigenvalues $\tau_i$ are
generally sorted in decreasing order so that $v_1^\prime y$ is the
`largest' principal component (that is, the principal component with
the largest variance). This is understandable when the goal is
expressed as capturing most of the total variation in the {\it first} few
principal components. But to avoid confusion with other
sections of this paper, we keep the convention of eigenvalues
being sorted in increasing order, and then express the goal as capturing
most of the total variation in the {\it largest} few principal components.

The $k$ largest
principal components in our notation are thus given by $v_p' y, \ldots, 
v_{p-k+1}' y$ $(k=1, \ldots, p)$. Their (cumulative) fraction captured of the
total variation contained in $y$, denoted by $f_k(\Sigma)$, is given by
\begin{equation}\label{eq:f-k-Sigma}
f_k(\Sigma) = \frac{\sum_{j=1}^{k} \tau_{p-j+1}}{\sum_{m=1}^p \tau_m}~,
\quad k = 1, \ldots, p~.
\end{equation} 

The most common rule in deciding how many principal components to
retain is to decide on a given fraction of the total
variation that one wants to capture, denoted by $f_{target}$, and to then retain the
largest $k$ principal components, where $k$ is the smallest integer
satisfying $f_k(\Sigma) \ge f_{target}$. Commonly chosen values of
$f_{target}$ are $70\%, 80\%, 90\%$, depending on the context. For
obvious reasons, this rule is known as the {\it
  cumulative-percentage-of-total-variation} rule.

There exist a host of other rules, either analytical or
graphical, such as Kaiser's rule or the scree plot; see
\citeauthor{jolliffe:2002} (2002, Section~6.1). The vast majority of these
rules are also solely based on the eigenvalues $(\tau_1, \ldots, \tau_p)$.

The problem is that generally the covariance matrix $\Sigma$ is
unknown. Thus, neither the (population) principal components $v_i^\prime y$ nor
their cumulative percentages of total variation $f_k(\Sigma)$ can be
used in practice.

The common solution is to replace 
$\Sigma$ with the sample covariance matrix $S$, computed from a random
sample $y_1, \ldots, y_n$, independent of $y$.
Let $((\lambda_{1},\ldots,\lambda_{p});(u_{1},\ldots,u_{p}))$ denote a system of
eigenvalues and eigenvectors of~$S$; it is assumed again that the
eigenvalues $\lambda_i$ are sorted in increasing order.
Then the (sample) principal components of $y$ are given
by $u_1^\prime y, \ldots, u_p^\prime y$. 
% The eigenvector $v_i$ is called {\it
%   vector of coefficients} or {\it loadings} for the $i$-th principal
% component ($i = 1, \ldots, p$).

The various rules in deciding how many (sample) principal components to retain
are now based on the sample
eigenvalues $\lambda_i$. For example, the
cumulative-percentage-of-total-variation rule retains the largest $k$
principal components, where $k$ is the smallest integer satisfying
$f_k(S) \ge f_{target}$, with
\begin{equation}\label{eq:f-k-S}
f_k(S) = \frac{\sum_{j=1}^{k} \lambda_{p-j+1}}{\sum_{m=1}^p \lambda_m}~,
\quad k = 1, \ldots, p~.
\end{equation} 

The pitfall in doing so, unless $p \ll n$, 
is that $\lambda_i$ is not good estimator of the variance of the
$i$th principal component. Indeed, the variance of the $i$th
principal component, $u_i^\prime y$,  is given by $u_i^\prime \Sigma u_i^{}$ rather than by
$\lambda_i = u_i^\prime S u_i^{}$. By design, for large values of $i$,
the estimator~$\lambda_i$ is upward biased for the true variance
$u_i^\prime \Sigma u_i^{}$, whereas for small values of $i$, it is
downward biased. In other words, the variances of the large principal
components are overestimated whereas the variances of the small
principal components are underestimated. The unfortunate consequence
is that most rules in deciding how many principal components to
retain, such as the cumulative-percentage-of-total-variation rule,
generally retain fewer principal components than really needed.

\subsection{Previous Approaches under Large-Dimensional Asymptotics}

All the previous approaches under large-dimensional asymptotics that
we are aware of impose some additional structure on the
estimation problem.

Most works assume a sparseness conditions on the eigenvectors $v_i$ or on
the covariance matrix $\Sigma$; see \cite{amini:2011} for a
comprehensive review.

\cite{mestre:2008b}, on the other hand and as discussed before, assumes that $\Sigma$ has only
$\bar M \ll p$ distinct
eigenvalues and further that the multiplicity of each of the $\bar M$
distinct eigenvalues is known (which implies that the number $\bar M$ is
known as well). Furthermore, he needs spectral separation.
In this restrictive setting, he is able to construct a consistent estimator of every
distinct eigenvalue and its associated eigenspace (that is, the
space spanned by all eigenvectors corresponding to a specific distinct
eigenvalue).

\subsection{Alternative Approach Based on Nonlinear Shrinkage}

Unlike previous approaches under large-dimensional asymptotics, 
we do not wish to impose additional
structure on the estimation problem. In such a setting, improved
estimators of the eigenvectors $v_i$ are not available and one must
indeed use the sample eigenvectors $u_i$ as loadings. Therefore, as in
common practice, the principal components used are the $u_i^\prime y$.

Ideally, the rules in deciding how many principal components to
retain should be based on the variances of the principal
components given by $u_i^\prime \Sigma u_i^{}$.
It is important to note that even if the
population eigenvalues $\tau_i$ were known, the rules should not be based on
them. This is because the population eigenvalues 
$\tau_i = v_i^\prime \Sigma v_i$ describe the variances of
the~$v_i^\prime y$, which are not available and thus not used.
It seems that this important point has not been realized
so far. Indeed, various authors have used PCA as a motivational example in the
estimation of the population eigenvalues $\tau_i$; for example, see
\cite{karoui:2008}, \cite{mestre:2008b}, and
\cite{yao:kammoun:najim:2012}. But unless the population eigenvectors
$v_i$ are available as well, using the $\tau_i$ is misleading.

Although, in the absence of additional structure, it is not possible to
construct improved principal components, it is possible to accurately
estimate the variances of the commonly-used principal components. This
is because the variance of the $i$th principal component is nothing
else than the finite-sample-optimal nonlinear shrinkage constant $d_i^*$; see
equation~\eqref{e:star}. Its oracle counterpart $d_i^{or}$ is given
in equation~\eqref{eq:oracle0} and the {\it bona fide} 
counterpart $\widehat d_i$ is given in equation~\eqref{eq:nonlin0-p-less-n}
in the case $p < n$ and in equation~\eqref{eq:nonlin0-p-greater-n} in
the case $p > n$.

Our solution then is to base the various rules in deciding how
many principal components to retain on the $\widehat d_i$ in place of
the unavailable $d_i^* = u_i^\prime \Sigma u_i^{}$. For example, the
cumulative-percentage-of-total-variation rule retains the $k$ largest
principal components, where $k$ is the smallest integer satisfying
$f_k(\widehat S) \ge f_{target}$, with
\begin{equation}\label{eq:f-k-S-hat}
f_k(\widehat S) = \frac{\sum_{j=1}^{k} \widehat d_{p-j+1}}{\sum_{m=1}^p
  \widehat d_m}~,
\quad k = 1, \ldots, p~.
\end{equation} 

\begin{remark}\label{rem:total-var}
\rm
We have taken the total variation to be $\sum_{m=1}^p d_m^*$,
and the variation attributable to the~$k$ largest principal components
to be $\sum_{j=1}^k d_{p-j+1}^*$. In general, the sample principal
components $u_i'y$ are not uncorrelated (unlike the population
principal components~$v_i' y$). This means that $u_i' \Sigma u_j^{}$ can be
non-zero for $i \neq j$. Nonetheless, even in this case, the variation attributable to
the $k$ largest principal components is still equal to $\sum_{j=1}^k
d_{p-j+1}^*$, as explained in Appendix~\ref{app:offdiagonal}.~\qed
\end{remark}

% It is noteworthy to point out that, under large-dimensional asymptotics,
% although it is not possible to accurately/consistently estimate the (entire) covariance
% matrix itself, in the absence of additional structure, it is
% nevertheless possibly to accurately/consistently estimate certain functionals of
% the covariance matrix, such as the variance of principal components.

% {\bf Not sure whether to leave in this last paragraph. On the one
%   hand, we can't really prove consistency; see below. On the other
%   hand, if we say ``accurately'' instead of ``consistently'', then
%   maybe that makes our estimator of the covariance matrix sound bad.}

While most applications of PCA seek the principal components with the
largest variances, there are also some applications of PCA that 
seek the principal components with the {\it smallest}
variances; see \citeauthor{jolliffe:2002} (2002, Section~3.4). In the
case $p > n$, a certain number of the $\lambda_i$ will be equal to
zero, falsely giving the impression that a certain number of the
smallest principal components have variance zero. Such applications
also highlight the use of replacing the $\lambda_i$ with our nonlinear
shrinkage constants $\widehat d_i$, which are always greater than
zero.

% {\bf Should we add a formal result here?
%   The problem is that we cannot really prove
%   that $f_k(\widehat S)$ is a consistent estimator of $f_k(\Sigma)$
%   since this corresponds one-to-one to proving that $\widehat S$ is
%   consistent for $S^*$ and all we can prove is that $\widehat S$ is
%   consistent for~$S^{or}$. So maybe only do proofs by simulation in the
%   context of PCA? Or has this other guy Pan been able to proof that
%   $\widehat S$ is also consistent for $S^*$ in the meantime?

%   What we can prove is that $f_k(\widehat S)$ is a consistent
%   estimator of $f_k(S^{or})$, a concept that would first have to be
%   introduced.

%   My vote would be for proofs by simulation only.}

\section{Monte Carlo Simulations}
\label{sec:monte-carlo}

In this section, we study the finite-sample performance of various
estimators in different settings.

\subsection{Estimation of Population Eigenvalues}
\label{ss:mc-H}

We first focus on estimating the eigenvalues of the population
covariance matrix, $\boldsymbol{\tau}_n$. Of major interest to us is
the case where all the eigenvalues are or can be distinct; but we also
consider the case where they are known or assumed to be grouped into a small number of high-multiplicity clusters.

\subsubsection{All Distinct Eigenvalues}
\label{subsub:continuous}

We consider the following estimators of $\boldsymbol{\tau}_n$.
\begin{packed_itemize} 
\item {\bf Sample:} The sample
  eigenvalues $\lambda_{n,i}$.
\item {\bf Lawley:} The bias-corrected
  sample eigenvalues using the formula of \citeauthor{lawley:1956}
  (1956, Section 4). This transformation may not be monotonic in finite samples. Therefore, we post-process it with an isotonic regression. 
\item {\bf El Karoui:} The estimator of \cite{karoui:2008}. It provides an estimator of $H_n$, not $\boldsymbol{\tau}_n$, 
so we derive estimates of the population eigenvalues using `smoothed'
quantiles in the spirit of equations
(\ref{eq:condexp})--(\ref{eq:inverse}).\footnote{We implemented this
  estimator to the best of our abilities, following the description
  in \cite{karoui:2008}. Despite several attempts, we
  were not able to obtain the original code.}
\item {\bf LW:} Our estimator $\widehat{\boldsymbol{\tau}}_n$ of Theorem~\ref{theo:estimator}.
\end{packed_itemize} 
It should be pointed out that the estimator of \cite{lawley:1956} is
designed to reduce the finite-sample bias of the sample eigenvalues
$\lambda_{n,i}$; it is not necessarily designed for consistent
estimation of $\boldsymbol{\tau}_n$ under large-dimensional asymptotics.

Let $\widetilde \tau_{n,i}$ denote a generic estimator of $\tau_{n,i}$.
The evaluation criterion is the dimension-normalized Euclidian distance between
 estimated eigenvalues $\widetilde{\boldsymbol{\tau}}_n$ and  population eigenvalues $\boldsymbol{\tau}_n$:
\begin{equation}\label{eq:criterion}
\frac{1}{p}\sum_{i=1}^p \left[\widetilde \tau_{n,i} - \tau_{n,i} \right]^2~,
\end{equation}
averaged over 1,000 Monte Carlo simulations in each
scenario. 

\medskip 
\noindent{\sc Convergence}

In the first design, the $i$th population eigenvalue is equal to
$\tau_{n,i}\defeq H^{-1} ((i-0.5)/p)$ $(i = 1, \ldots, p)$, where $H$ is
given by the distribution of $1 + 10 W$,
and $W \sim \mbox{Beta}(1, 10)$; this distribution is right-skewed
and resembles in shape an exponential distribution. 
The distribution of the random
variates comprising the $n \times p$ data matrix $X_n$ is real Gaussian.
We fix the concentration at $p/n = 0.5$ and vary the dimension from $p=30$ to
$p=1,000$. The results are displayed in Figure~\ref{fig:conv1}. 
\begin{center}
\captionsetup{type=figure}
\includegraphics[scale=\scale]{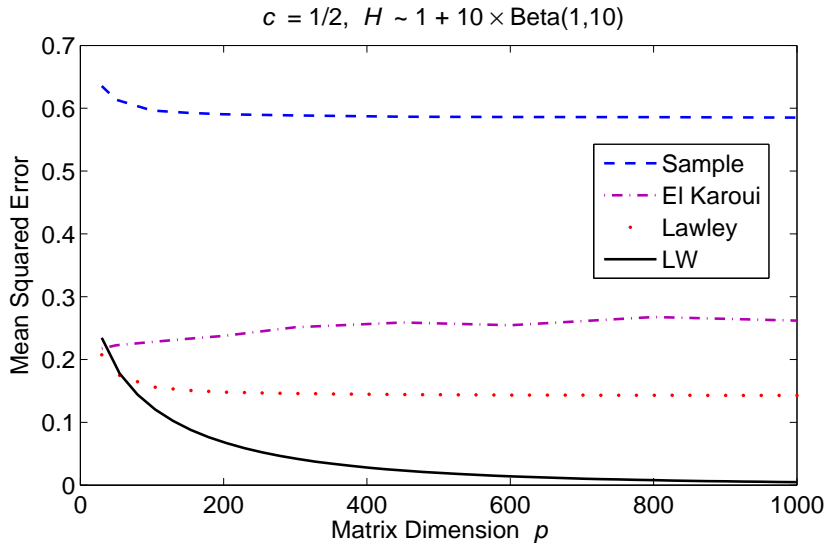}
\captionof{figure}{Convergence of estimated eigenvalues to population
  eigenvalues in the case where the sample covariance matrix is
  nonsingular.}\label{fig:conv1}
\end{center}
% software references for Olivier [do not delete this comment]:
% convergence_mse02.eps, optim19, Results05, msedistance02.m
% tau51_*.m, nek04_*.m, 29 Dec 2012
It can be seen that
the empirical mean squared error for LW converges to zero,
which is in agreement with the proven consistency of
Theorem~\ref{theo:estimator}. For all the other estimators, the
average distance from $\boldsymbol{\tau}_n$ appears bounded away from zero.
% in particular, the failure of El Karoui to converge to zero appears on
% contradiction to the theoretical results provided by
% \cite{karoui:2008}.
This simulation also shows that dividing by~$p$ is indeed the
appropriate normalization for the Euclidian norm in equation
(\ref{eq:individual}), as it drives a wedge between estimators such as
the sample eigenvalues that are not consistent and
$\widehat{\boldsymbol{\tau}}_n$, which is consistent.

The second design is similar to the first design, except that we fix
the concentration at $p/n = 2$ and now vary the sample size from
$n=30$ to $n=1,000$ instead of the dimension. In this design, the
sample covariance matrix is always singular. The results are displayed
in Figure~\ref{fig:conv2} and are qualitatively similar. Again, LW
is the only estimator
that appears to be consistent. Notice the vertical scale: the
difference between El Karoui and LW  is of
the same order of magnitude as in Figure \ref{fig:conv1}, but Sample
and Lawley are vastly more erroneous now.

\begin{center}
\captionsetup{type=figure}
\includegraphics[scale=\scale]{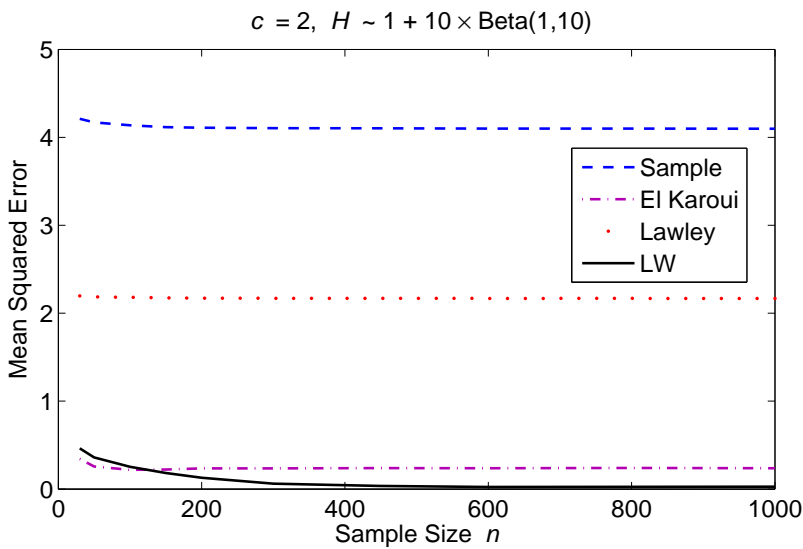}
\captionof{figure}{Convergence of estimated eigenvalues to population eigenvalues in the case where the sample covariance matrix is singular.}
\label{fig:conv2}
\end{center}

% Software versions for Olivier's reference: 7 November 2012, mp\_07, optim\_21, Results06, tau53\_*.m, nek04.m, msedistance03.m, convergence\_mse03.eps

\medskip
\noindent{\sc Condition Number}

In the third design, the focus is on the condition number. 
The $i$th population eigenvalue is still 
$\tau_{n,i}\defeq H^{-1} ((i-0.5)/p)$ $(i = 1, \ldots, p)$, but 
$H$ is now given by the distribution of $a + 10W$, where $W \sim
\mbox{Beta}(1,10)$, and $a \in [0, 7]$. As a result, 
the smallest eigenvalue approaches $a$, and the previously-used
distribution for $H$ is included as a special case when $a =1$. The condition
number decreases in $a$ from approximately $10,000$ to $2.4$. 

We use
$n=1,600$ and $p=800$, so that $p/n = 0.5$. The random
variates are still real Gaussian. The results are displayed in Figure~\ref{fig:cond}. It can
be seen that Sample and Lawley perform quite well for
values of $a$ near zero (that is, for very large condition numbers)
but their performance gets worse as $a$ increases (that is, as the
condition number decreases). On the other hand, the performance of El
Karoui is more stable across all values of $a$, though relatively bad. The
performance of LW is uniformly the best and also stable across $a$.

\begin{center}
\captionsetup{type=figure}
\includegraphics[scale=\scale]{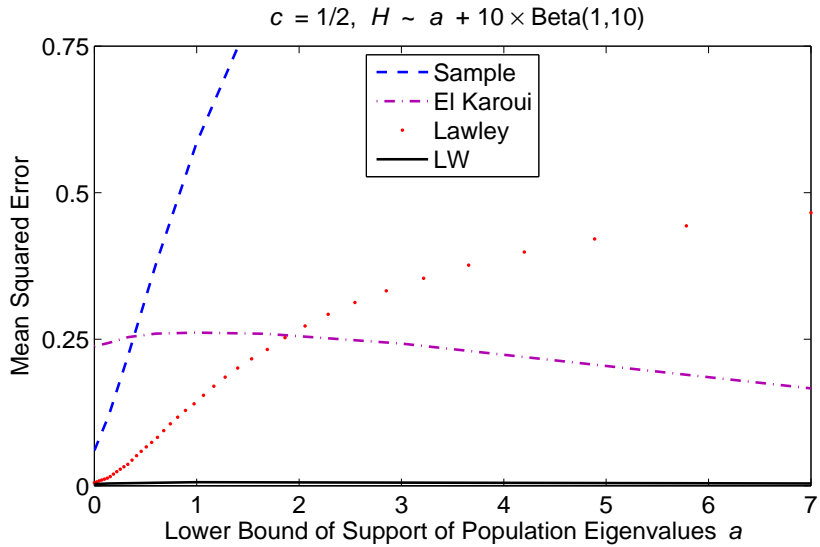}
\captionof{figure}{Effect of the condition number on the mean squared error between estimated and population eigenvalues.}
\label{fig:cond}
\end{center}

% Software versions for Olivier's reference: mp\_07, optim32\Results20, optim32\tau63_*.m
% optim_22\nek07_*.m, optim_22\msedistance06.m, optim22\condition11.eps

\medskip
\noindent{\sc Shape of the Distribution}

In the fourth design, we consider a wide variety of shapes of $H$,
which is now given by the distribution of $1 + 10W$, where $W$ follows
a Beta distribution with parameters
$\{(1,1),(1,2),(2,1),$
$(1.5,1.5),(0.5,0.5),(5,5),(5,2),(2,5)\}$; see Figure~7 of
\cite{ledoit:wolf:2012} for a graphical representation of the
corresponding densities. Always again, the $i$th population eigenvalue is 
$\tau_{n,i}\defeq H^{-1} ((i-0.5)/p)$ $(i = 1, \ldots, p)$.

We use $n=1,600$ and $p=800$, so that $p/n = 0.5$. The random
variates are real Gaussian. The results are presented in Table~\ref{table:shape}.
It can be seen that LW is uniformly best and Sample is uniformly worst. There is no clear-cut
ranking for the remaining two estimators. On average, Lawley is
second best, followed by El Karoui.

\begin{center}
\captionsetup{type=table}
$$\begin{array}{|c||c|c|c|c|c|}\hline
\mbox{Parameters} &  \mbox{LW} &  \mbox{Sample} &  \mbox{El Karoui} &  \mbox{Lawley} \\ \hline\hline
(1,1) & 0.15 & 6.70 & 2.65 & 0.66 \\ \hline
(1,2) & 0.06 & 2.58 & 1.65 & 0.27 \\ \hline
(2,1) & 0.16 & 15.59 & 2.23 & 2.61 \\ \hline
(1.5,1.5) & 0.09 & 7.07 & 2.03 & 0.93 \\ \hline
(0.5,0.5) & 0.08 & 7.04 & 2.87 & 0.53 \\ \hline
(5,5) & 0.08 & 9.52 & 1.02 & 2.13 \\ \hline
(5,2) & 0.12 & 20.93 & 1.39 & 4.90 \\ \hline
(2,5) & 0.08 & 2.59 & 0.87 & 0.46 \\ \hline \hline
\mbox{Average} & 0.10 & 9.00 & 1.84 & 1.56 \\ \hline
\end{array}
$$
\captionof{table}{Mean squared error between estimated
  and population eigenvalues.}
\label{table:shape}
\end{center}

% Software versions for Olivier's reference: mp\_07, optim\_23, Results08, tau57\_*.m, nek08\_*.m, msedistance09.m

\medskip
\noindent{\sc Heavy Tails}

So far, the variates making up the data matrix $X_n$ always had a
Gaussian distribution. It is also of interest to consider a
heavy-tailed distribution instead. We return to the first design with
$n = 1600$ and $p=800$, so that $p/n = 0.5$. In addition to
the Gaussian distribution, which can be viewed as a $t$-distribution
with infinite degrees of freedom, we also consider a the $t$-distribution
with three degrees of freedom (scaled to have unit variance).
The results are presented in Table~\ref{table:fattail}. 
It can be seen that all estimators perform worse when
the degrees of freedom are changed from infinity to three, but 
LW is still by far the best.
\begin{center}
\captionsetup{type=table}
$$\begin{array}{|c||c|c|c|c|c|}\hline
\mbox{Degrees of Freedom} & \mbox{LW} & \mbox{Sample} & \mbox{El Karoui} & \mbox{Lawley} \\ \hline\hline
3 & 0.21 & 4.97 & 4.02 & 4.41 \\ \hline
\infty & 0.01 & 0.59 & 0.27 & 0.14 \\ \hline
\end{array}
$$
\captionof{table}{Mean squared error between
estimated and population eigenvalues.}
\label{table:fattail}
\end{center}

\subsubsection{Clustered Eigenvalues}

We are mainly interested in the case where the population eigenvalues are or can be distinct, but it is also worthwhile seeing how (an adapted version of) our
estimator of $\boldsymbol{\tau}_n$ compares to the one of \cite{mestre:2008b} in the 
setting where the population eigenvalues are known or assumed to be grouped into a small number of high-multiplicity clusters.

Let $\gamma_1<\gamma_2<\dots<\gamma_{\bar{M}}$ denote the set of
pairwise different eigenvalues of the population covariance matrix
$\Sigma$, where $\bar{M}$ is the number of distinct population
eigenvalues ($1\leq\bar{M}<p$). Each of the eigenvalues $\gamma_j$ has
known multiplicity $K_j \; (j=1,\ldots,\bar{M})$, so that
$p=\sum_{j=1}^{\bar{M}}K_j$. (Knowing the multiplicities of the
eigenvalues $\gamma_j$ comes from knowing their masses $m_j$ in the
limiting spectral distribution $H$, as assumed in \cite{mestre:2008b}:
 $K_j / p = m_j$.)

Then the optimization problem in Theorem 2.1 becomes:
\begin{align}
&\widehat{\boldsymbol{\gamma}}_n \defeq \argmin_{(\gamma_1,\gamma_2,\dots,\gamma_{\bar{M}})\in[0,\infty)^{\bar{M}}}\;\frac{1}{p}
\sum_{i=1}^p\left[\lambda_{n,i}-q_{n,p}^i(\mathbf{t})\right]^2
\label{eq:mod1} \\
\mbox{subject to:}\quad\mathbf{t}&=(
\underbrace{\gamma_1,\ldots,\gamma_1}_{K_1\,\mbox{times}},
\underbrace{\gamma_2,\ldots,\gamma_2}_{K_2\,\mbox{times}},\ldots,
\underbrace{\gamma_{\bar{M}},\ldots,\gamma_{\bar{M}}}_{K_{\bar{M}}\,\mbox{times}}
)' \label{eq:mod2} \\
& \! \gamma_1 < \gamma_2 < \ldots < \gamma_{\bar M} \label{eq:mod3}
\end{align}
We consider the following estimators of $\boldsymbol{\tau}_n$.

\begin{packed_itemize} 
\item {\bf Traditional:} $\gamma_j$ is estimated by the average
  of all corresponding sample \mbox{eigenvalues~$\lambda_{n,i}$}; under the condition of
  spectral separation assumed in \cite{mestre:2008b}, it is known
  which~$\gamma_j$ corresponds to which $\lambda_{n,i}$.
\item {\bf Mestre:} The estimator defined in
  \citeauthor{mestre:2008b} (2008, Theorem~3).
\item {\bf LW:} Our modified estimator as defined in \eqref{eq:mod1}--\eqref{eq:mod3}.
\end{packed_itemize} 
The mean squared error criterion \eqref{eq:criterion} specializes in this setting to
$$
\sum_{j=1}^{\bar M} m_j \, (\widehat \gamma_j - \gamma_j)^2~.
$$
We report the average MSE over 1,000 Monte Carlo simulations in
each scenario. 

\medskip
\noindent{\sc Convergence}

The first design is based on Tables I and II of \cite{mestre:2008b}. The
distinct population eigenvalues are $(\gamma_1, \gamma_2, \gamma_3, \gamma_4) = (1,
  7, 15, 25)$ with respective multiplicities $(K_1, K_2, K_3, K_4) = (p/2, p/4, p/8, p/8)$.
The distribution of the random
variates comprising the $n \times p$ data matrix $X_n$ is circular
symmetric complex Gaussian, as in \cite{mestre:2008b}.
We fix the concentration at $p/n = 0.32$ and vary the dimension from $p=8$ to
$p=1,000$; the lower end $p = 8$ corresponds to Table~I in
\cite{mestre:2008b}, while the upper end $p=1,000$ corresponds to
Table~II in \cite{mestre:2008b}.
The results are displayed in Figure~\ref{fig:discrete}. It can be seen that
the average MSE of both Mestre and LW converges to zero,
and that the performance of the two estimators is
nearly indistinguishable. On the other hand, Traditional is seen to be
inconsistent, as its MSE remains bounded away from zero.

\begin{center}
\captionsetup{type=figure}
\includegraphics[scale=\scale]{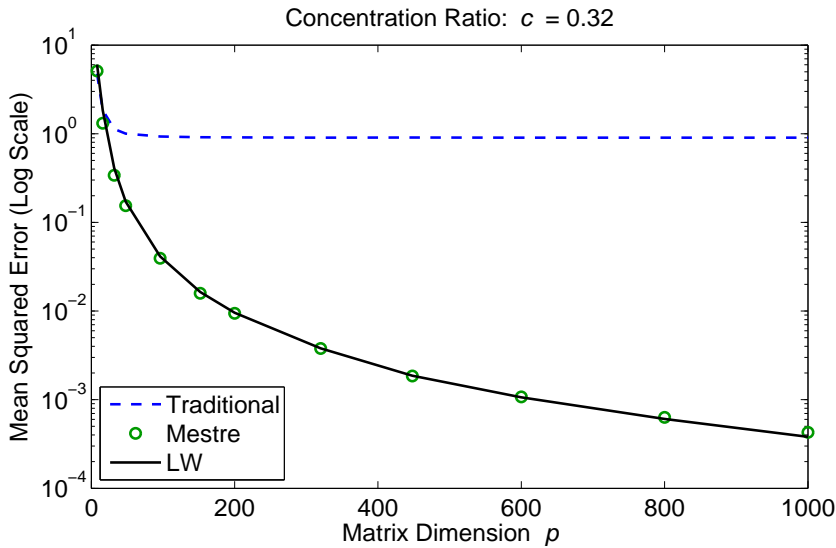}
\captionof{figure}{Convergence of estimated eigenvalues to population eigenvalues when eigenvalues are grouped into a small number of high-multiplicity clusters.}
\label{fig:discrete}
\end{center}

% Software versions for Olivier's reference: 30 December 2012, mp\_07, optim\_24, Results12, mestre09.m, plotter17.m, mestre17.eps

\medskip
\noindent{\sc Performance When One Eigenvalue Is Isolated}

The second design is based on Table III of \cite{mestre:2008b}. 
The
distinct population eigenvalues are $(\gamma_1, \gamma_2, \gamma_3, \gamma_4) = (1,
  7, 15, 25)$ with multiplicities $(K_1, K_2, K_3, K_4) = (160, 80, 79, 1)$.
There is a single `isolated' large
eigenvalue. 
The distribution of the random
variates comprising the $n \times p$ data matrix $X_n$ is circular
symmetric complex Gaussian, as in \cite{mestre:2008b}.
We use $n = 1,000$ and $p = 320$, so that $p/n = 0.32$. The averages and the standard deviations of the estimates $\widehat
\gamma_j$ over 10,000 Monte Carlos simulations
are presented in Table~\ref{table:3}; note here that the numbers for
Traditional and Mestre have been directly copied from Table~III of \cite{mestre:2008b}.
The inconsistency of Traditional is again apparent. In terms of
estimating $(\gamma_1, \gamma_2, \gamma_3)$, the performance of Mestre
and LW is nearly indistinguishable. In terms of estimating $\gamma_4$,
Mestre has a smaller bias (in absolute value) while LW has a smaller
standard deviation; combining the two criteria yields a mean squared
error of $(25 - 24.9892)^2 + 1.0713^2 = 1.1478$ for Mestre and a mean
squared error of $(25 - 24.9238)^2 + 0.8898^2 = 0.7976$ for LW.

\begin{center}
\captionsetup{type=table}
$$\begin{array}{|c|c||c|c||c|c||c|c|}\cline{3-8}
\multicolumn{2}{c}{}&\multicolumn{2}{|c||}{\mbox{Traditional}}&\multicolumn{2}{|c||}{\mbox{Mestre}}&\multicolumn{2}{|c|}{\mbox{LW}}\\  \hline 
\mbox{Eigenvalue} & \mbox{Multiplicity} & \mbox{Mean} & \mbox{Std.~Dev.} & \mbox{Mean} & \mbox{Std.~Dev.} & \mbox{Mean} & \mbox{Std.~Dev.}\\ \hline \hline 
\gamma_1=1  & 160 &  0.8210 & 0.0023 &  0.9997 & 0.0032 &  1.0006 & 0.0034 \\ \hline
\gamma_2=7  &  80 &  6.1400 & 0.0208 &  6.9942 & 0.0343 &  7.0003 & 0.0319 \\ \hline
\gamma_3=15 &  79 & 16.1835 & 0.0514 & 14.9956 & 0.0681 & 14.9995 & 0.0580 \\ \hline
\gamma_4=25 &   1 & 28.9104 & 0.7110 & 24.9892 & 1.0713 & 24.9238 & 0.8898 \\ \hline 
\end{array}
$$ 
\captionof{table}{Empirical mean and standard
  deviation of the eigenvalue estimator of
  \protect\cite{mestre:2008b}, sample eigenvalues, and the proposed
  estimator. The first six columns are copied from \mbox{Table~III} of
  \protect\cite{mestre:2008b}. Results are based on $10,000$
  Monte Carlo simulations with circularly symmetric complex Gaussian
  random variates. 
% Mean Squared Error is computed as an average across Monte Carlo simulations: $\sum_{i=1}^{10,000}\left[(\widehat{\lambda}_1^i-1)^2\times160+
% (\widehat{\lambda}_2^i-7)^2\times80+(\widehat{\lambda}_3^i-15)^2\times79+(\widehat{\lambda}_4^i-25)^2\times1\right]/(320\times10,000)$.
}
\label{table:3}
\end{center}

%Software versions for Olivier's reference: 03 September 2012, mp\_07, optim\_24, Results11, mestre08.m

\subsection{Covariance Matrix Estimation}
\label{ss:mc-sigma}

As detailed in
Section~\ref{ss:oracle}, the finite-sample optimal estimator in the
class of rotation-equiva\-riant estimators is given by~$S_n^*$ as
defined in~\eqref{e:optimal-fs}. As the benchmark, we use 
the linear shrinkage estimator of
\cite{ledoit:wolf:2004a} instead of the sample
covariance matrix. We do this because the linear shrinkage estimator
has become the {\em de facto} standard
among leading researchers because of its simplicity, accuracy, and good
conditioning properties. It has been used in several fields of
statistics, such as linear regression with a large number of regressors
\citep{Anatolyev2012368}, linear discriminant analysis
\citep{PedroDuarteSilva20112975}, factor analysis \citep{Lin2012448},
unit root tests \citep{Demetrescu201210}, and vector autoregressive
models \citep{Huang2011}, among others. Beyond pure statistics, the linear shrinkage
estimator has been applied in finance for portfolio selection
\citep{Tsagaris20121651} and tests of asset pricing models
\citep{Khan200855}; in signal processing for cellular phone
transmission \citep{Nguyen2011176} and radar detection
\citep{Wei2011}; and in biology for neuroimaging
\citep{NIPS2010_1054}, genetics \citep{Lin2012631}, cancer research
\citep{Pyeon20074605}, and psychiatry \citep{Markon2010273}. It has
also been used in such varied applications as physics \citep{Pirkl2012533}, chemistry
\citep{Guo20123880}, climatology \citep{Ribes2009707}, oil exploration \citep{Saetrom2012152}, road safety research \citep{Haufe2011}, 
 etc. In summary, the
comparatively poor performance of the sample covariance matrix and the
popularity of the linear shrinkage estimator justify taking the latter
as the benchmark.

The improvement of the nonlinear shrinkage estimator $\widehat{S}_n$ over
the linear shrinkage estimator of \cite{ledoit:wolf:2004a}, denoted by
$\overline S_n$, will be measured by how closely this
estimator approximates the finite-sample optimal estimator 
$S^*_n$ relative to $\overline S_n$.
 More specifically, we report the Percentage Relative
Improvement in Average Loss (PRIAL), which is defined as
\begin{equation}\label{e:prial}
\mbox{PRIAL} \defeq {\rm
  PRIAL}(\widehat \Sigma_n) \defeq 100\times
\left\{1-\frac{\E \Bigl [\bigl \|\widehat{\Sigma}_n-S^*_n\bigr \|_F^2\Bigr ]}
{\E\left[\bigl \|\overline S_n-S^*_n \big \|_F^2\right]}\right\} \%~,
\end{equation}
where $\widehat \Sigma_n$ is an arbitrary estimator of $\Sigma_n$.
By definition, the PRIAL of $\overline S_n$ is 0\% while the PRIAL of $S_n^*$ is 100\%.

We consider the following estimators of $\Sigma_n$.
\begin{packed_itemize} 
\item {\bf LW (2012) Estimator:} The nonlinear shrinkage estimator of
  \cite{ledoit:wolf:2012}; this version only works for the case $p < n$.
\item {\bf New Nonlinear Shrinkage Estimator:} The new nonlinear
  shrinkage estimator of Section~\ref{ss:nonlin}; this version works across both cases $p < n$
  and $p > n$.
\item{\bf Oracle:} The (infeasible) oracle estimator of Section~\ref{ss:oracle}.
\end{packed_itemize} 

\medskip
\noindent{\sc Convergence}

In our design, 20\% of the population eigenvalues are equal to 1,
40\% are equal to 3, and 40\% are equal to 10. This is a particularly
interesting and difficult example introduced and analyzed in detail by
\cite{bai:silverstein:1998}; it has also been used in previous Monte Carlo
simulations by \cite{ledoit:wolf:2012}.
The distribution of the random
variates comprising the $n \times p$ data matrix $X_n$ is real Gaussian.
We study convergence of the various estimators by keeping the
concentration $p/n$ fixed while increasing the sample size $n$.
We consider the three cases $p/n = 0.5, 1, 2$; as discussed in
Remark~\ref{rem:c=1}, the case $p/n = 1$ is not covered by the mathematical
treatment. 
The results are displayed in Figure~\ref{fig:prial},
which shows empirical PRIAL's across 1,000 Monte
Carlo simulations (one panel for each case $p/n = 0.5, 1, 2$).
It can be seen that the new nonlinear shrinkage
estimator always outperforms linear shrinkage with its PRIAL
converging to 100\%, though slower than the oracle estimator. As
expected, the relative improvement over the linear shrinkage estimator
is inversely related to the concentration ratio; also see Figure~4 of
\cite{ledoit:wolf:2012}. In the case $p < n$, it can also be seen that
the new nonlinear shrinkage estimator slightly outperforms the
earlier nonlinear shrinkage estimator of
\cite{ledoit:wolf:2012}. Last but not least, although the case $p=n$
is not covered by the mathematical treatment, it is also dealt with
successfully in practice by the new nonlinear shrinkage estimator.

\begin{center}
\captionsetup{type=figure}
\includegraphics[scale=\scale]{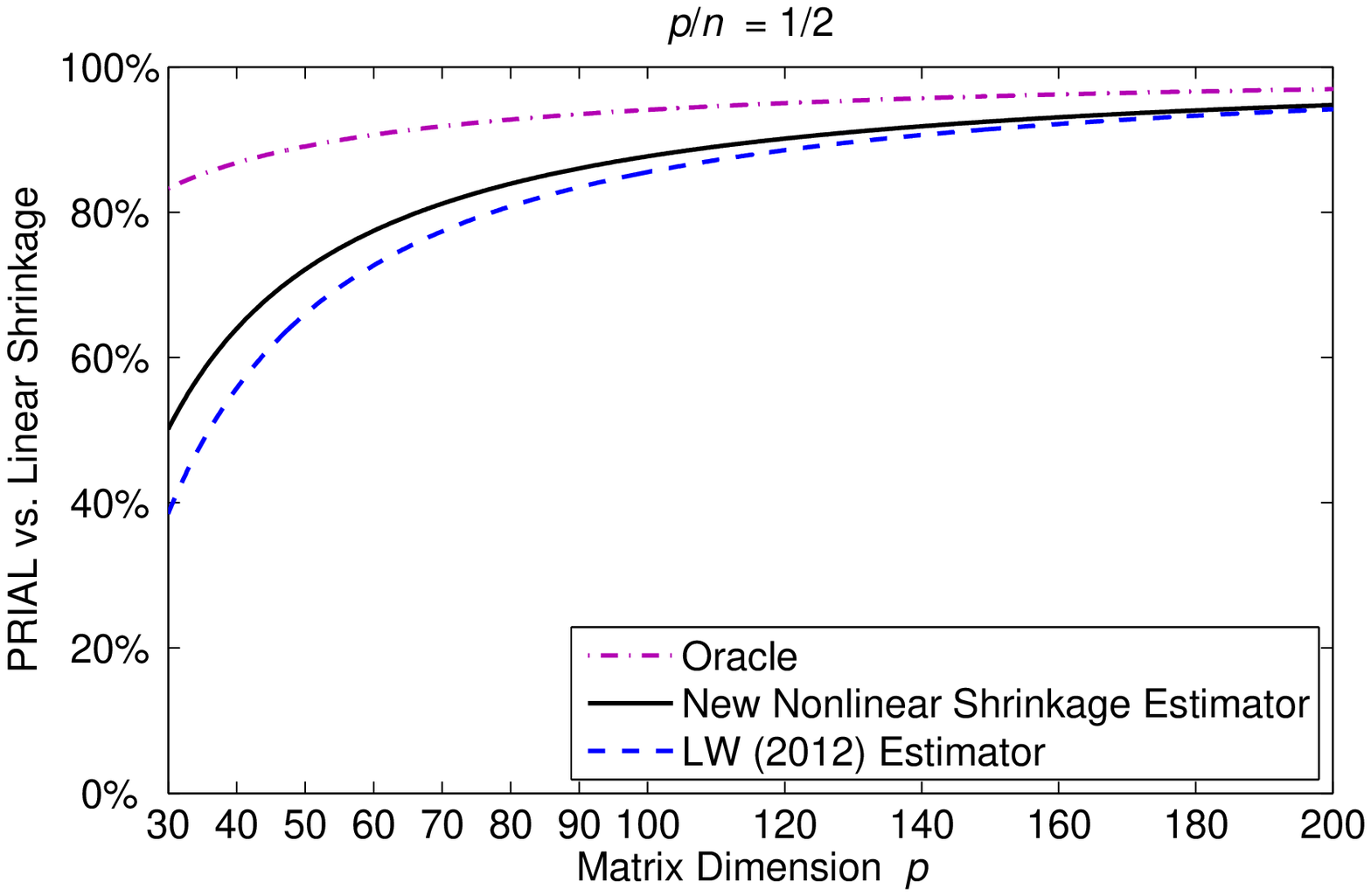}
% Software versions for Olivier's reference: 22 November 2012, convergence10, results, caller3\_fun1.m, mp\_07, convergence\_plotter7.m, convergence\_7.eps
\includegraphics[scale=\scale]{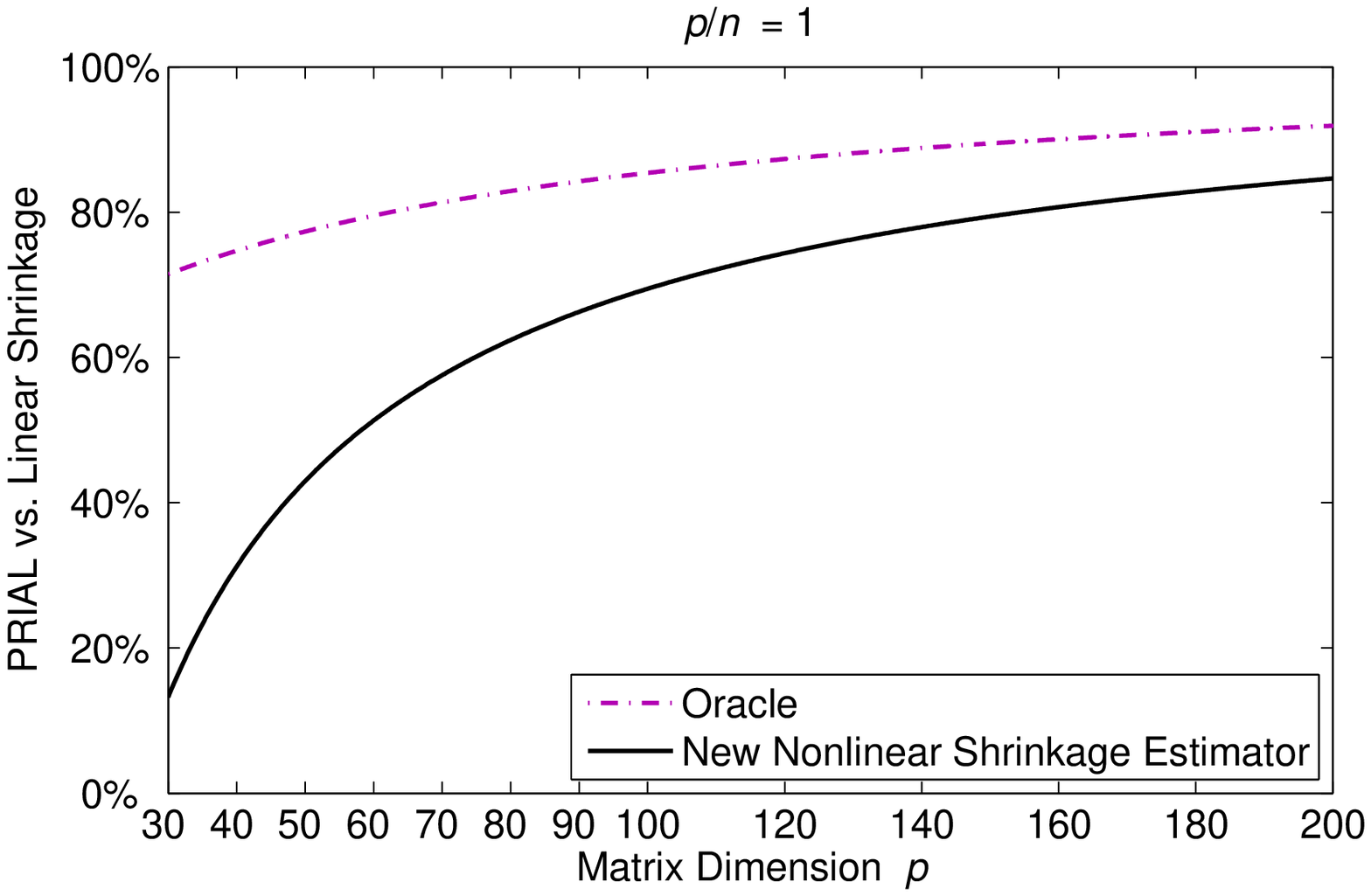}
% Software versions for Olivier's reference: 25 Nov.~2012, optim\_27, Results15, tau59\_*.m, mp\_08, convergence\_plotter11.m, p\_equals\_n_v2.eps
\includegraphics[scale=\scale]{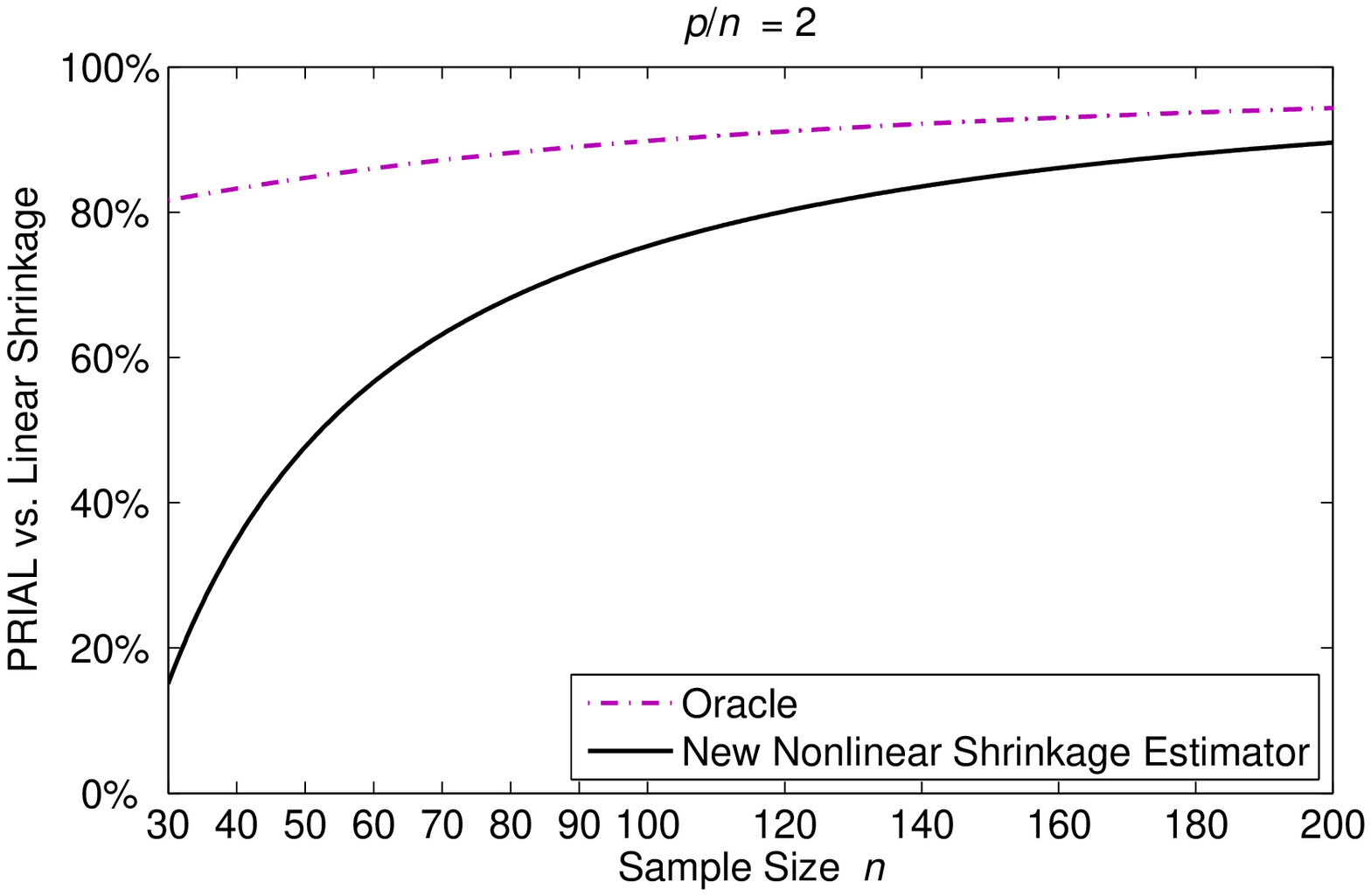}
% Software versions for Olivier's reference: 22 Nov.~2012, optim\_26, Results14, tau58\_*.m, mp\_07, convergence\_plotter16.m, prial16.eps
\captionof{figure}{Percentage Improvement in Average Loss (PRIAL) according to the Frobenius norm of nonlinear versus linear shrinkage estimation of the covariance matrix.}
\label{fig:prial}
\end{center}

\subsection{Principal Component Analysis}
\label{ss:mc-pca}

Recall that in Section \ref{sec:pca} on principal component analysis we dropped the first subscript $n$ always, and so the same will be done in this section.

In our design, the $i$th population eigenvalue is equal to $\tau_{i}=H^{-1} ((i-0.5)/p)$ $(i = 1, \ldots, p)$, where $H$ is given by the distribution of $1 + 10 W$,
and $W \sim \mbox{Beta}(1, 10)$
The distribution of the random
variates comprising the $n \times p$ data matrix $X_n$ is Gaussian.
We consider the two cases $(n = 200, p=100)$ and $(n=100, p = 200)$,
so the concentration is $p/n = 0.5$ or $p/n =2$.

Let $y \in \R^p$ be a random vector with covariance matrix $\Sigma$, drawn independently from the sample covariance matrix $S$.
The out-of-sample variance of the $i$th (sample) principal component, $u_i'y$, is given by
$d_i^* \defeq u_i'\Sigma u_i^{}$; see \eqref{e:star}.
By our convention, the $d_i^*$ are sorted in increasing order.

We consider the following estimators of $d_i^*$.
\begin{packed_itemize} 
\item {\bf Sample:} The estimator of $d_i^*$ is the $i$th sample
  eigenvalue, $\lambda_i$.
\item {\bf Population:} The estimator of $d_i^*$ is the $i$th
  population eigenvalue, $\tau_i$; this estimator is not feasible but is
  included for educational purposes nevertheless.
\item {\bf LW:} The estimator of $d_i^*$ is the nonlinear
  shrinkage quantity $\widehat d_i$ as given in equation~\eqref{eq:nonlin0-p-less-n}
in the case $p < n$ and in equation~\eqref{eq:nonlin0-p-greater-n} in
the case $p > n$.
\end{packed_itemize} 
Let $\widetilde d_i$ be a generic estimator of $d_i^*$. 
First, we are plotting
$$
\widetilde f_k \defeq \frac{\sum_{j=1}^{k} \widetilde d_{p-j+1}}{\sum_{m=1}^p
  \widetilde d_m}
$$
as a function of $k$, averaged over 1,000 Monte Carlo simulations.
The quantity $\widetilde f_k$  serves as an
estimator of $f_k$, the fraction of the total variation in~$y$ that is explained by
the $k$ largest principal components:
$$
f_k \defeq \frac{\sum_{j=1}^{k} d_{p-j+1}^*}{\sum_{m=1}^p d_m^*}
$$
The results are displayed in Figure~\ref{fig:pca} (one panel for each case
 $p/n = 0.5, 2$.)
The upward bias of Sample is  apparent,
while LW is very close to the Truth. Moreover, Population is also
upward biased (though not as much as Sample): the important message is
that even if the population eigenvalues were known, they should not be
used to judge the variances of the (sample) principal components. As
expected, the differences between Sample and LW increase with the
concentration $p/n$; the same is true for the differences between
Population and~LW.

\begin{center}
\captionsetup{type=figure}
\includegraphics[scale=\scale]{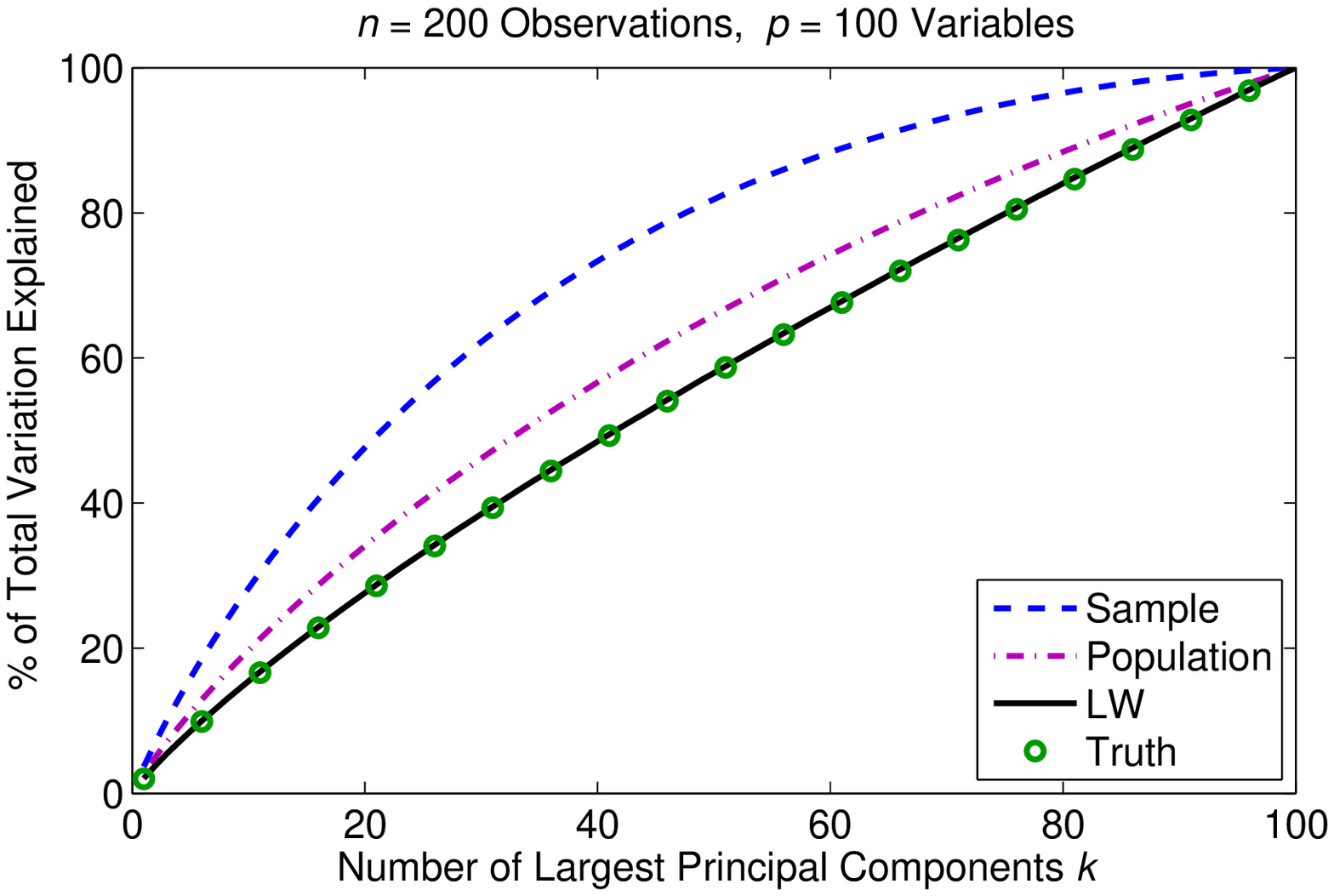}
% Software versions for Olivier's reference: 09 November 2012, optim\_19, tau51\_100\_200.m, mp\_07, pcaplotter07.m, pca07.eps
\captionsetup{type=figure}
\includegraphics[scale=\scale]{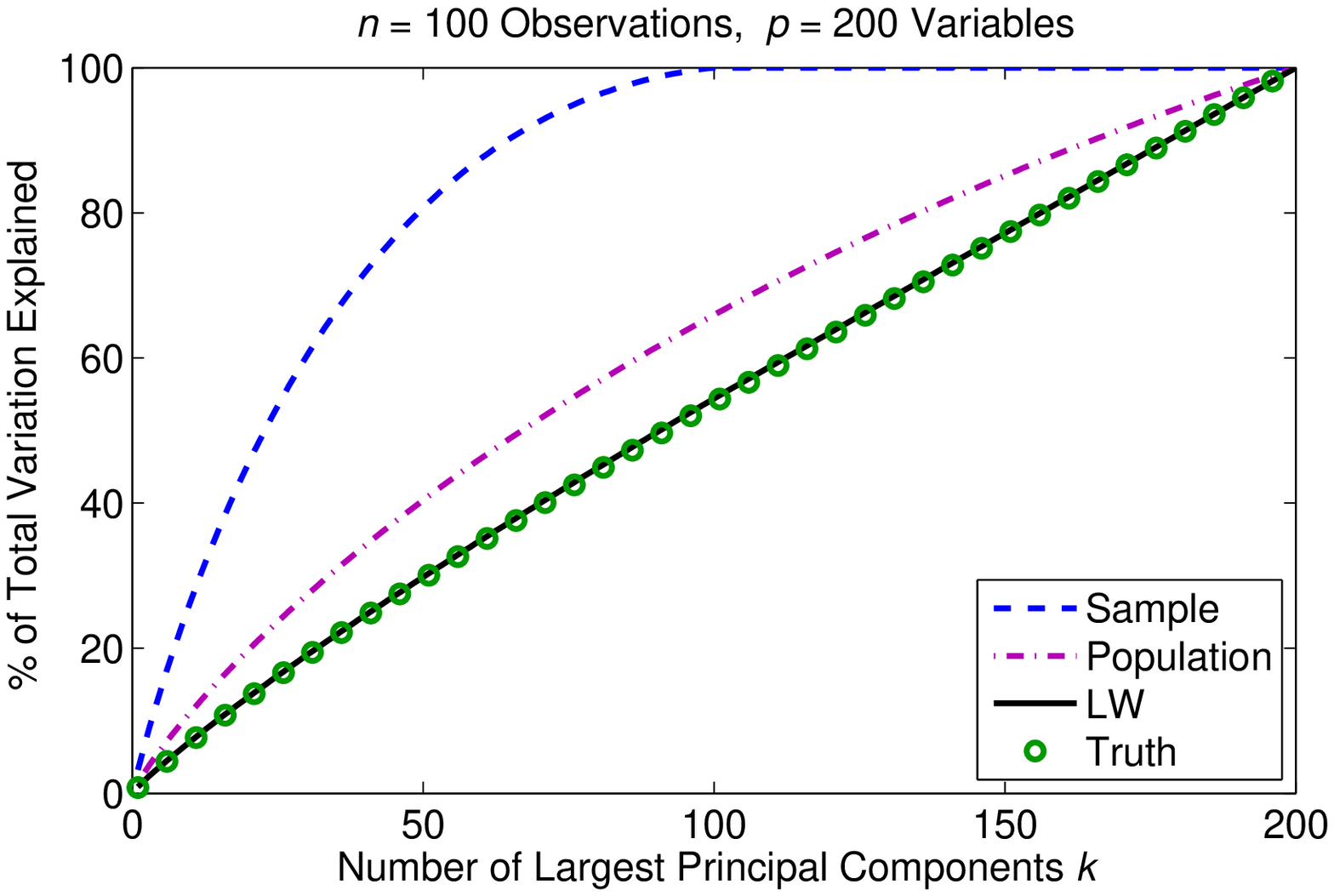}
% Software versions for Olivier's reference: 09 November 2012, optim\_19, tau51\_200\_100.m, mp\_07, pcaplotter08.m, pca08.eps
\captionof{figure}{Comparison between different estimators of the percentage of total variation explained by the top principal components.}
\label{fig:pca}
\end{center}

Figure~\ref{fig:pca} shows how close the estimator $\widetilde f_k$
is to the truth $f_k$ on average. But it does not necessarily answer how
close the cumulative-percentage-of-total-variation rule based on
$\widetilde f_k$ is to the rule based on $f_k$. For a given percentage
$(q \times 100) \%$, with $q \in (0, 1)$, let
$$
k(q) \defeq \min\{k : f_k \ge q\} \qquad \mbox{and} \qquad
\widetilde k(q) \defeq \min\{k : \widetilde f_k \ge q\}~.
$$
In words, $k(q)$ is the (smallest) number of the largest principal
components that must be retained to explain $(q \times
100)\%$ of the total variation and $\widetilde k(q)$ is an estimator
of this quantity.

We are then also interested in the Root Mean Squared Error (RMSE)
of $\widetilde k(q)$, defined as
$$
\mbox{RMSE} \defeq \sqrt{\E \left [\bigl (\widetilde k(q) - k(q) \bigr )^2 \right ]}~,
$$
for the values of $q$ most commonly used in practice, namely
$q = 0.7, 0.8, 0.9$. We compute empirical RMSEs across 1,000 Monte Carlo
simulations. The results are presented in Table~\ref{table:pca-rmse}. It
can be seen that in each scenario, Sample has the largest RMSE and LW
has the smallest RMSE; in particular LW is highly accurate not only
relative to Sample but also in an absolute sense.
(The quantity $k(q)$ is a random variable, since it depends on the
sample eigenvalues~$u_i$; but the general magnitude of $k(q)$ for the
various scenarios can be judged from Figure~\ref{fig:pca})

\begin{center}
\captionsetup{type=table}
$$\begin{array}{|c||c|c|c|c|c||c|c|c|}\hline
\mbox{$q$} & \mbox{Sample} & \mbox{Population} & \mbox{LW} & & 
\mbox{$q$} & \mbox{Sample} & \mbox{Population} & \mbox{LW} \\ \hline\hline
\multicolumn{4}{|c|}{\rule[-2mm]{0mm}{6.5mm} {\it n}  = 200, {\it p} =100} &
\multicolumn{1}{c|}{\rule[-2mm]{0mm}{6.5mm} } &
\multicolumn{4}{c|}{\rule[-2mm]{0mm}{6.5mm} {\it n} = 100, {\it p} = 200} 
\\ \hline
70\% & 26.9 & 9.1 & 0.8 & & 70\% & 96.5 & 25.4 & 1.4    \\ \hline
80\% & 27.9 & 8.0 & 0.8 & & 80\% & 107.3 & 21.0 & 0.9    \\ \hline
90\% & 24.4 & 5.0 & 0.6 & & 90\% & 114.0 & 13.0 & 0.7    \\ \hline
\end{array}
$$
\captionof{table}{Empirical Root Mean Squared Error (RMSE) of various
  estimates, $\widetilde k(q)$, of the
  number of largest principal components that must be retained, $k(q)$, to
  explain $(q \times 100)\%$ of the total variation. Based on
  1,000 Monte Carlo simulations.}
\label{table:pca-rmse}
\end{center}

\section{Empirical Application}
\label{sec:emp-app}

As an empirical application, we study principal component analysis
(PCA) in the context of stock return data. 
Principal components of a return vector of a cross section of $p$
stocks are used
for risk analysis and portfolio selection by finance practitioners;
for example, see \cite{roll:ross:1980} and
\cite{connor:korajczyk:1993}.

We use the $p = 30, 60,
240, 480$ largest stocks, as measured by their market value at the
beginning of 2011, that have a complete return history from January 2002 until
 December 2011. As is customary in many financial applications, such
 as portfolio selection, we use monthly data. Consequently, the
 sample size for the ten-year history is $n = 120$ and the
 concentration is $p/n = 0.25, 0.5, 2, 4$.

It is of crucial interest how much of the total variation in the
$p$-dimensional return vector is explained by the $k$ largest
principal components. We compare the approach based on the sample
covariance matrix, defined in equation~\eqref{eq:f-k-S} and denoted by
Sample, to that of nonlinear shrinkage, defined in
equation~\ref{eq:f-k-S-hat} and denoted by LW. The results are
displayed in Figure~\ref{fig:pca-stocks}. It can be seen that Sample is
overly optimistic compared to LW and, as expected, the differences
between the two methods increase with the concentration $p/n$.

\begin{center}
\captionsetup{type=figure}
\includegraphics[width=6in,height=6in]{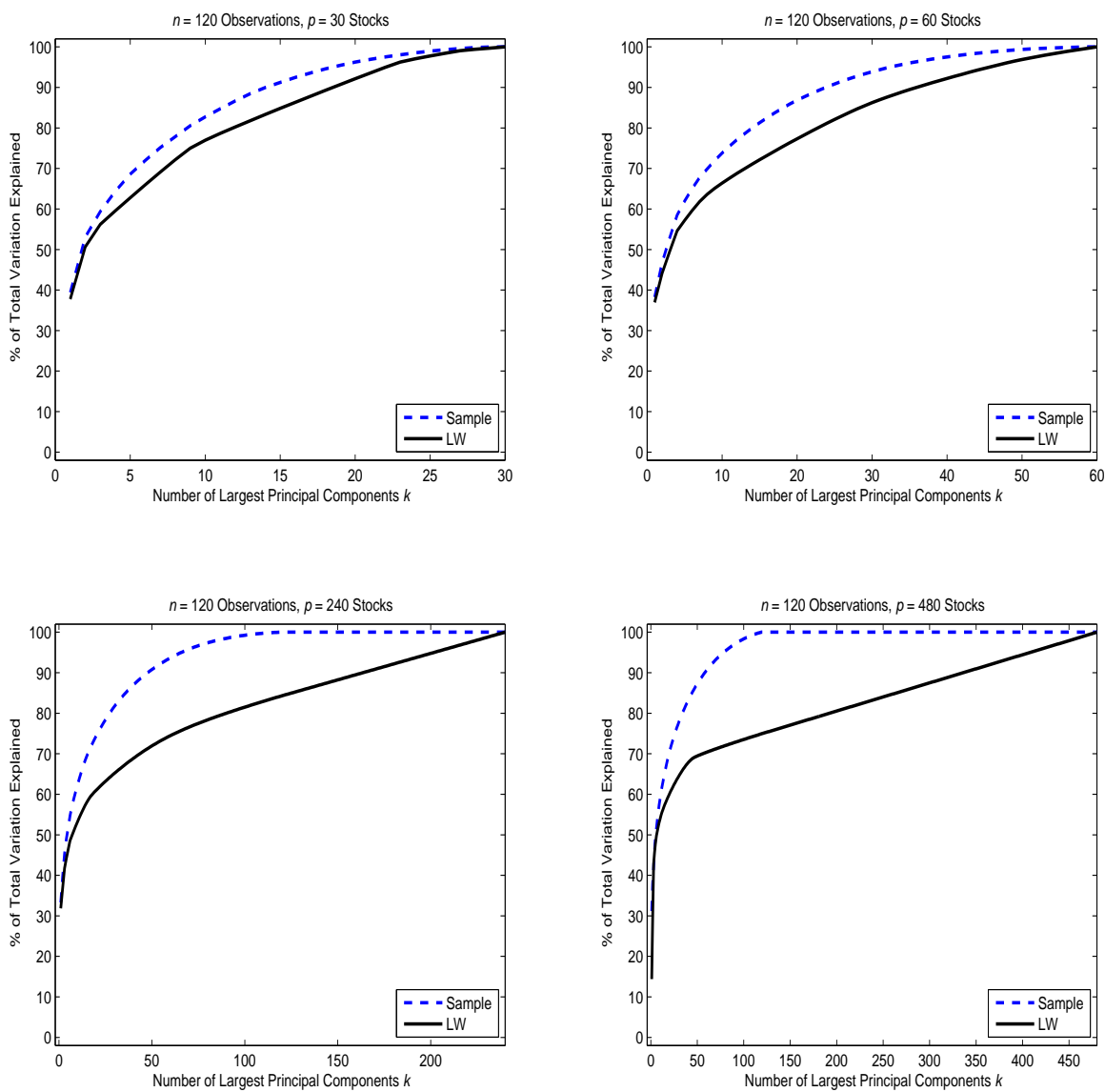}
\captionof{figure}{Percentage of total variation explained by the $k$ largest
  principal components of stock returns: estimates based on the sample covariance
  matrix (Sample) compared to those based on nonlinear shrinkage
  (LW).}
\label{fig:pca-stocks}
\end{center}

In addition to the visual analysis, the differences can also be
presented via the cumulative-percentage-of-total-variation rule to
decide how many of the largest principal components to retain; see
Section~\ref{sec:pca}. The results are presented in
Table~\ref{table:pca}. It can be seen again that Sample is overly
optimistic compared to LW and retains much fewer principal
components. Again, as expected, the
differences between the two methods increase with the concentration $p/n$.

\begin{center}
\captionsetup{type=table}
$$\begin{array}{|c||c|c|c|c||c|c|}\hline
\mbox{$f_{Target}$} & \mbox{LW} & \mbox{Sample} & &
\mbox{$f_{Target}$} & \mbox{LW} & \mbox{Sample} \\ \hline\hline
\multicolumn{3}{|c|}{\rule[-2mm]{0mm}{6.5mm} {\it n}  = 120, {\it p} = 30} &
\multicolumn{1}{c|}{\rule[-2mm]{0mm}{6.5mm} } &
\multicolumn{3}{c|}{\rule[-2mm]{0mm}{6.5mm} {\it n} = 120, {\it p} = 60} 
\\ \hline
70\% & 8 & 6  & & 70\% & 14 & 9     \\ \hline
80\% & 12 & 9 & & 80\% & 23 & 15    \\ \hline
90\% & 19 & 15& & 90\% & 36 & 24    \\ \hline \hline
\multicolumn{3}{|c|}{\rule[-2mm]{0mm}{6.5mm} {\it n} = 120, {\it p} = 240} &
\multicolumn{1}{c|}{\rule[-2mm]{0mm}{6.5mm} } &
\multicolumn{3}{c|}{\rule[-2mm]{0mm}{6.5mm} {\it n} = 120, {\it p} = 480} 
\\ \hline
70\% & 44 & 16& & 70\% & 56 & 20    \\ \hline
80\% & 91 & 28& & 80\% & 193 & 34    \\ \hline
90\% & 164 & 48& & 90\% & 337 & 59   \\ \hline 
\end{array}
$$
\captionof{table}{Number of $k$ largest principal components to retain according to the
  cumulative-percentage-of-total-variation rule. The rule is based
  either on the sample covariance matrix (Sample) or on nonlinear
  shrinkage (LW).}
\label{table:pca}
\end{center}

\section{Conclusion}
\label{sec:conclusion}

The analysis of large-dimensional data sets is becoming more and more
common. For many statistical problems, the classic textbook
methods no longer work well in such settings. 
Two cases in point are covariance matrix estimation and principal
component analysis, both cornerstones of multivariate analysis.

The classic estimator of the covariance matrix is the sample
covariance matrix. It is unbiased and the maximum-likelihood estimator
under normality. But when the dimension is not small compared to the
sample size, the sample covariance matrix
 contains too much estimation error and is ill-conditioned; when the dimension
is larger than the sample size, it is not even invertible anymore.

The variances of the principal components (which are obtained from the
sample eigenvectors) are no longer accurately estimated by the sample
eigenvalues. In particular, the sample eigenvalues overestimate the
variances of the `large' principal components. As a result the common
rules in determining how many principal components to retain generally
select too few of them.

In the absence of strong structural assumption on the true covariance
matrix, such as sparseness or a factor model, a
common remedy for both statistical problems is {\it nonlinear
  shrinkage} of the sample eigenvalues. The optimal shrinkage formula
delivers an estimator of the covariance matrix
that is finite-sample optimal with respect to the Frobenius norm 
in the class of rotation-equivariant
estimators. The {\it same} shrinkage
formula also gives the variances of the (sample) principal
components. It is noteworthy that the optimal shrinkage formula is
different from the population eigenvalues: even if they were
available, one should not use them for the ends of covariance matrix
estimation and PCA.

Unsurprisingly, the optimal nonlinear shrinkage formula is not
available, since it depends on population quantities. But an
asymptotic counterpart, denoted {\it oracle shrinkage}, can be
estimated consistently. In this way, {\it bona fide} nonlinear
shrinkage estimation of covariance matrices and improved PCA result. 

The key to the consistent estimation of the oracle shrinkage is the
consistent estimation of the population eigenvalues. This
problem is challenging and interesting in its own right, solving a host
of additional statistical problems. Our proposal to this end is
the first one that does not make strong assumptions on the
distribution of the population eigenvalues, has proven consistency, and also works well in practice.

Extensive Monte Carlo simulations have established that our methods
have desirable finite-sample properties and outperform methods that
have been previously suggested in the literature.

\bibliographystyle{apalike}
\bibliography{../wolf}

\begin{thebibliography}{}

\bibitem[Amini, 2011]{amini:2011}
Amini, A.~A. (2011).
\newblock High-dimensional principal component analysis.
\newblock Technical Report {UCB/EECS}-2011-104, Department of Electrical
  Engineering and Computer Sciences, University of California at Berkeley.

\bibitem[Anatolyev, 2012]{Anatolyev2012368}
Anatolyev, S. (2012).
\newblock Inference in regression models with many regressors.
\newblock {\em Journal of Econometrics}, 170(2):368--382.

\bibitem[Bai and Silverstein, 1998]{bai:silverstein:1998}
Bai, Z.~D. and Silverstein, J.~W. (1998).
\newblock No eigenvalues outside the suppport of the limiting spectral
  distribution of large-dimensional random matrices.
\newblock {\em Annals of Probability}, 26(1):\mbox{316--345}.

\bibitem[Bai and Silverstein, 1999]{bai:silverstein:1999}
Bai, Z.~D. and Silverstein, J.~W. (1999).
\newblock Exact separation of eigenvalues of large-dimensional sample
  covariance matrices.
\newblock {\em Annals of Probability}, 27(3):1536--1555.

\bibitem[Bai and Silverstein, 2010]{bai:silverstein:2010}
Bai, Z.~D. and Silverstein, J.~W. (2010).
\newblock {\em Spectral Analysis of Large-Dimensional Random Matrices}.
\newblock Springer, New York, second edition.

\bibitem[Bickel and Freedman, 1981]{bickel:freedman:1981}
Bickel, P.~J. and Freedman, D.~A. (1981).
\newblock Asymptotic theory for the bootstrap.
\newblock {\em Annals of Statistics}, 9(6):1196--1217.

\bibitem[Bickel and Levina, 2008]{bickel:levina:2008}
Bickel, P.~J. and Levina, E. (2008).
\newblock Regularized estimation of large covariance matrices.
\newblock {\em Annals of Statistics}, 36(1):199--227.

\bibitem[Connor and Korajczyk, 1993]{connor:korajczyk:1993}
Connor, G. and Korajczyk, R.~A. (1993).
\newblock A test for the number of factors in an approximate factor model.
\newblock {\em Journal of Finance}, 48:1263--1291.

\bibitem[Demetrescu and Hanck, 2012]{Demetrescu201210}
Demetrescu, M. and Hanck, C. (2012).
\newblock A simple nonstationary-volatility robust panel unit root test.
\newblock {\em Economics Letters}, 117(1):10--13.

\bibitem[{El Karoui}, 2008]{karoui:2008}
{El Karoui}, N. (2008).
\newblock Spectrum estimation for large dimensional covariance matrices using
  random matrix theory.
\newblock {\em Annals of Statistics}, 36(6):2757--2790.

\bibitem[Fan et~al., 2008]{fan:fan:lv:2008}
Fan, J., Fan, Y., and Lv, J. (2008).
\newblock High dimensional covariance matrix estimation using a factor model.
\newblock {\em Journal of Econometrics}, 147(1):186--197.

\bibitem[Gill et~al., 2002]{gill:murray:saunders:2002}
Gill, P.~E., Murray, W., and Saunders, M.~A. (2002).
\newblock {SNOPT}: An {SQP} algorithm for large-scale constrained optimization.
\newblock {\em SIAM Journal on Optimization}, 12(4):979--1006.

\bibitem[Guo et~al., 2012]{Guo20123880}
Guo, S.-M., He, J., Monnier, N., Sun, G., Wohland, T., and Bathe, M. (2012).
\newblock Bayesian approach to the analysis of fluorescence correlation
  spectroscopy data {II}: Application to simulated and in vitro data.
\newblock {\em Analytical Chemistry}, 84(9):3880--3888.

\bibitem[Haufe et~al., 2011]{Haufe2011}
Haufe, S., Treder, M., Gugler, M., Sagebaum, M., Curio, G., and Blankertz, B.
  (2011).
\newblock {EEG} potentials predict upcoming emergency brakings during simulated
  driving.
\newblock {\em Journal of Neural Engineering}, 8(5).

\bibitem[Hotelling, 1933]{hotelling:1933}
Hotelling, H. (1933).
\newblock Analysis of a complex of statistical variables into principal
  components.
\newblock {\em Journal of Educational Psychology}, 24(6):417--441, 498--520.

\bibitem[Huang and Schneider, 2011]{Huang2011}
Huang, T.-K. and Schneider, J. (2011).
\newblock Learning auto-regressive models from sequence and non-sequence data.
\newblock In Shawe-Taylor, J., Zemel, R., Bartlett, P., Pereira, F., and
  Weinberger, K., editors, {\em Advances in Neural Information Processing
  Systems 24}, pages 1548--1556. The {MIT} Press, Cambridge.

\bibitem[Jolliffe, 2002]{jolliffe:2002}
Jolliffe, I.~T. (2002).
\newblock {\em Principal Component Analysis}.
\newblock Springer, New York, second edition.

\bibitem[Khan, 2008]{Khan200855}
Khan, M. (2008).
\newblock Are accruals mispriced? {E}vidence from tests of an intertemporal
  capital asset pricing model.
\newblock {\em Journal of Accounting and Economics}, 45(1):55--77.

\bibitem[Lawley, 1956]{lawley:1956}
Lawley, D.~N. (1956).
\newblock A general method for approximating to the distribution of likelihood
  ratio criteria.
\newblock {\em Biometrika}, 43(3/4):295--303.

\bibitem[Ledoit and P\'ech\'e, 2011]{ledoit:peche:2011}
Ledoit, O. and P\'ech\'e, S. (2011).
\newblock Eigenvectors of some large sample covariance matrix ensembles.
\newblock {\em Probability Theory and Related Fields}, 150(1--2):233--264.

\bibitem[Ledoit and Wolf, 2004]{ledoit:wolf:2004a}
Ledoit, O. and Wolf, M. (2004).
\newblock A well-conditioned estimator for large-dimensional covariance
  matrices.
\newblock {\em Journal of Multivariate Analysis}, 88(2):365--411.

\bibitem[Ledoit and Wolf, 2012]{ledoit:wolf:2012}
Ledoit, O. and Wolf, M. (2012).
\newblock Nonlinear shrinkage estimation of large-dimensional covariance
  matrices.
\newblock {\em Annals of Statistics}, 40(2):1024--1060.

\bibitem[Lin and Bentler, 2012]{Lin2012448}
Lin, J. and Bentler, P. (2012).
\newblock A third moment adjusted test statistic for small sample factor
  analysis.
\newblock {\em Multivariate Behavioral Research}, 47(3):448--462.

\bibitem[Lin et~al., 2012]{Lin2012631}
Lin, J.-A., Zhu, H.~b., Knickmeyer, R., Styner, M., Gilmore, J., and Ibrahim,
  J. (2012).
\newblock Projection regression models for multivariate imaging phenotype.
\newblock {\em Genetic Epidemiology}, 36(6):631--641.

\bibitem[Mar{\v{c}}enko and Pastur, 1967]{marcenko:pastur:1967}
Mar{\v{c}}enko, V.~A. and Pastur, L.~A. (1967).
\newblock Distribution of eigenvalues for some sets of random matrices.
\newblock {\em Sbornik: Mathematics}, 1(4):457--483.

\bibitem[Markon, 2010]{Markon2010273}
Markon, K. (2010).
\newblock Modeling psychopathology structure: A symptom-level analysis of axis
  {I} and {II} disorders.
\newblock {\em Psychological Medicine}, 40(2):273--288.

\bibitem[Mestre, 2008]{mestre:2008b}
Mestre, X. (2008).
\newblock Improved estimation of eigenvalues and eigenvectors of covariance
  matrices using their sample estimates.
\newblock {\em {IEEE} Transactions on Information Theory}, 54(11):5113--5129.

\bibitem[Nguyen et~al., 2011]{Nguyen2011176}
Nguyen, L., Rheinschmitt, R., Wild, T., and Brink, S. (2011).
\newblock Limits of channel estimation and signal combining for multipoint
  cellular radio (comp).
\newblock In {\em Wireless Communication Systems (ISWCS), 2011 8th
  International Symposium on}, pages 176--180. IEEE.

\bibitem[Pearson, 1901]{pearson:1901}
Pearson, K. (1901).
\newblock On line and planes of closest fit to systems of points in space.
\newblock {\em Philosophical Magazine (Series 6)}, 2(11):559--572.

\bibitem[Pedro Duarte~Silva, 2011]{PedroDuarteSilva20112975}
Pedro Duarte~Silva, A. (2011).
\newblock Two-group classification with high-dimensional correlated data: A
  factor model approach.
\newblock {\em Computational Statistics and Data Analysis}, 55(11):2975--2990.

\bibitem[Perlman, 2007]{perlman:2007}
Perlman, M.~D. (2007).
\newblock {\em STAT 542: Multivariate Statistical Analysis}.
\newblock University of Washington (On-Line Class Notes), Seattle, Washington.

\bibitem[Pirkl et~al., 2012]{Pirkl2012533}
Pirkl, R., Remley, K., and Patan√©, C. (2012).
\newblock Reverberation chamber measurement correlation.
\newblock {\em IEEE Transactions on Electromagnetic Compatibility},
  54(3):533--545.

\bibitem[Pyeon et~al., 2007]{Pyeon20074605}
Pyeon, D., Newton, M., Lambert, P., Den~Boon, J., Sengupta, S., Marsit, C.,
  Woodworth, C., Connor, J., Haugen, T., Smith, E., Kelsey, K., Turek, L., and
  Ahlquist, P. (2007).
\newblock Fundamental differences in cell cycle deregulation in human
  papillomavirus-positive and human papillomavirus-negative head/neck and
  cervical cancers.
\newblock {\em Cancer Research}, 67(10):4605--4619.

\bibitem[Rajaratnam et~al., 2008]{bala:et:al:2008}
Rajaratnam, B., Massam, H., and Carvalho, C.~M. (2008).
\newblock Flexible covariance estimation in graphical {Gaussian} models.
\newblock {\em Annals of Statistics}, 36(6):2818--2849.

\bibitem[Ribes et~al., 2009]{Ribes2009707}
Ribes, A., Aza√≠s, J.-M., and Planton, S. (2009).
\newblock Adaptation of the optimal fingerprint method for climate change
  detection using a well-conditioned covariance matrix estimate.
\newblock {\em Climate Dynamics}, 33(5):707--722.

\bibitem[Roll and Ross, 1980]{roll:ross:1980}
Roll, R. and Ross, S.~A. (1980).
\newblock An empirical investigation of the arbitrage pricing theory.
\newblock {\em Journal of Finance}, 35:1073--1103.

\bibitem[S{\ae}trom et~al., 2012]{Saetrom2012152}
S{\ae}trom, J., Hove, J., Skjervheim, J.-A., and Vab{\o}, J. (2012).
\newblock Improved uncertainty quantification in the ensemble {K}alman filter
  using statistical model-selection techniques.
\newblock {\em SPE Journal}, 17(1):152--162.

\bibitem[Silverstein, 1995]{silverstein:1995}
Silverstein, J.~W. (1995).
\newblock Strong convergence of the empirical distribution of eigenvalues of
  large-dimensional random matrices.
\newblock {\em Journal of Multivariate Analysis}, 55:331--339.

\bibitem[Silverstein and Bai, 1995]{silverstein:bai:1995}
Silverstein, J.~W. and Bai, Z.~D. (1995).
\newblock On the empirical distribution of eigenvalues of a class of
  large-dimensional random matrices.
\newblock {\em Journal of Multivariate Analysis}, 54:175--192.

\bibitem[Silverstein and Choi, 1995]{silverstein:choi:1995}
Silverstein, J.~W. and Choi, S.~I. (1995).
\newblock Analysis of the limiting spectral distribution of large-dimensional
  random matrices.
\newblock {\em Journal of Multivariate Analysis}, 54:295--309.

\bibitem[Stein, 1975]{stein:1975}
Stein, C. (1975).
\newblock Estimation of a covariance matrix.
\newblock {Rietz lecture}, 39th Annual Meeting IMS. Atlanta, Georgia.

\bibitem[Stein, 1986]{stein:1986}
Stein, C. (1986).
\newblock Lectures on the theory of estimation of many parameters.
\newblock {\em Journal of Mathematical Sciences}, 34(1):1373--1403.

\bibitem[Tsagaris et~al., 2012]{Tsagaris20121651}
Tsagaris, T., Jasra, A., and Adams, N. (2012).
\newblock Robust and adaptive algorithms for online portfolio selection.
\newblock {\em Quantitative Finance}, 12(11):1651--1662.

\bibitem[Varoquaux et~al., 2010]{NIPS2010_1054}
Varoquaux, G., Gramfort, A., Poline, J.-B., and Thirion, B. (2010).
\newblock Brain covariance selection: better individual functional connectivity
  models using population prior.
\newblock In Lafferty, J., Williams, C. K.~I., Shawe-Taylor, J., Zemel, R., and
  Culotta, A., editors, {\em Advances in Neural Information Processing Systems
  23}, pages 2334--2342. The {MIT} Press, Cambridge.

\bibitem[Wei et~al., 2011]{Wei2011}
Wei, Z., Huang, J., and Hui, Y. (2011).
\newblock Adaptive-beamforming-based multiple targets signal separation.
\newblock In {\em Signal Processing, Communications and Computing (ICSPCC),
  2011 IEEE International Conference on}, pages 1--4. IEEE.

\bibitem[Yao et~al., 2012]{yao:kammoun:najim:2012}
Yao, J., Kammoun, A., and Najim, J. (2012).
\newblock Estimation of the covariance mat\-rix of large dimensional data.
\newblock Technical report, T\'{e}l\'{e}com ParisTech.
\newblock Available online at {\tt http://arxiv.org/abs/1201.4672}.

\end{thebibliography}

\newpage

\begin{appendix}

\renewcommand{\baselinestretch}{1.09}

\section{Mathematical Proofs}

\begin{lemma}\label{lem:grid-con}
Let $\{G_n\}$ and $G$ be c.d.f.'s on the real line and assume that there
exists a compact interval that contains the support of $G$ as well as the
support of $G_n$ for all $n$ large enough. For $0 < \alpha < 1$, let $G_n^{-1} (\alpha)$
denote an $\alpha$ quantile of $G_n$ and let $G^{-1}(\alpha)$ denote
an $\alpha$ quantile of $G$. Let $\{K_n\}$ be a sequence of integers with $K_n \to
\infty$; further, let $t_{n,k} \defeq (k-0.5)/K_n$ \mbox{$(k = 1, \ldots, K_n)$}.

\smallskip
Then $G_n \Rightarrow G$ if and only if
\begin{equation} \label{eq:conv-sum-finite}
\frac{1}{K_n} \sum_{k=1}^{K_n} \bigl [G_n^{-1} (t_{n,k}) - G^{-1}
(t_{n,k}) \bigr ]^2 \to 0~.
\end{equation} 
\end{lemma}

{\sc Proof. } First, since there exists a compact interval that contains the
support of $G$ as well as the support of $G_n$ for all $n$ large
enough, weak convergence of $G_n$ to $G$ implies that also the
second moment of $G_n$ converges to the second moment of $G$. 

Second, we claim that, under the given set of assumptions, the convergence
\eqref{eq:conv-sum-finite} is equivalent to the convergence
\begin{equation} \label{eq:conv-int}
\int_0^1 \bigl [G_n^{-1} (t) - G^{-1} (t) \bigr ]^2 dt \to 0~.
\end{equation} 
If this claim is true, the proof of the lemma
then follows from Lemmas~8.2 and~8.3(a) of
\cite{bickel:freedman:1981}. 
We are, therefore, left to show the claim that the 
convergence~\eqref{eq:conv-sum-finite} is equivalent to the convergence
\eqref{eq:conv-int}. 

We begin by showing that the 
convergence~\eqref{eq:conv-sum-finite} implies the convergence
\eqref{eq:conv-int}. One can make the following decomposition:
\begin{align}
\int_0^1 \bigl [G_n^{-1} (t) - G^{-1} (t) \bigr ]^2 dt = & 
\; \; \; \; \; \int_0^{t_{n,1}} \bigl [G_n^{-1} (t) - G^{-1} (t) \bigr
]^2 dt \nonumber \\
& + \int_{t_{n,1}}^{t_{n,K_n}} \bigl [G_n^{-1} (t) - G^{-1} (t)
\bigr]^2 dt \nonumber \\
& + \int_{t_{n,K_n}}^{1} \bigl [G_n^{-1} (t) - G^{-1} (t) \bigr ]^2 dt~. \nonumber
\end{align}
Let $C$ denote the length of the compact interval that contains the
support of $G$ and $G_n$ (for all~$n$ large enough). Also, note that
$t_{n,1}-0 = 1 - t_{n,K_n} = 0.5 / K_n$.
Then, for all $n$
large enough, we can use the trivial bound
$$
\int_0^{t_{n,1}} \bigl [G_n^{-1} (t) - G^{-1} (t) \bigr]^2 dt
+  \int_{t_{n,K_n}}^{1} \bigl [G_n^{-1} (t) - G^{-1} (t) \bigr ]^2 dt
\le  \frac{0.5}{K_n} C^2 + \frac{0.5}{K_n} C^2 = \frac{C^2}{K_n}~.
$$
Combining this bound with the previous decomposition results in
$$
\int_0^1 \bigl [G_n^{-1} (t) - G^{-1} (t) \bigr ]^2 dt \le 
\frac{C^2}{K_n} + \int_{t_{n,1}}^{t_{n,K_n}} \bigl [G_n^{-1} (t) - G^{-1} (t)
\bigr]^2 dt~,
$$
and we are left to show that
$$
\int_{t_{n,1}}^{t_{n,K_n}} \bigl [G_n^{-1} (t) - G^{-1} (t)
\bigr]^2 dt \to 0~.
$$

For any $k = 1, \ldots, K_n -1$, noting that $t_{n,k+1} - t_{n,k} = 1/K_n$,
\begin{align}
\int_{t_{n,k}}^{t_{n,k+1}} \bigl [G_n^{-1} (t) - G^{-1} (t) \bigr ]^2 dt
    & \le  
\frac{1}{K_n} \sup_{t_{n,k} \le t \le t_{n,k+1}} \bigl [G_n^{-1} (t) -
G^{-1} (t) \bigr ]^2 \nonumber \\
   & \le \frac{\bigl [G_n^{-1} (t_{n,k+1}) - G^{-1} (t_{n,k}) \bigr ]^2
       + \bigl [G_n^{-1} (t_{n,k}) - G^{-1} (t_{n,k+1}) \bigr ]^2 }{K_n}~,  \nonumber
\end{align}
where the last inequality follows from the fact that both $G_n^{-1}$ and
$G^{-1}$ are (weakly) increasing functions. As a result,
\begin{align}
\int_{t_{n,1}}^{t_{n,K_n}} \bigl [G_n^{-1} (t) - G^{-1} (t) \bigr ]^2
dt 
= & \sum_{k=1}^{K_n-1} \int_{t_{n,k}}^{t_{n,k+1}} \bigl [G_n^{-1} (t)
- G^{-1} (t) \bigr ]^2 \nonumber \\
\le & 
\; \; \; \; \; \frac{1}{K_n} \sum_{k=1}^{K_n-1} \bigl [G_n^{-1} (t_{n,k+1}) -
G^{-1}(t_{n,k}) \bigr ]^2 \label{eq:sum1} \\
& +
\frac{1}{K_n} \sum_{k=1}^{K_n-1} \bigl [G_n^{-1} (t_{n,k}) -
G^{-1}(t_{n,k+1}) \bigr ]^2 \label{eq:sum2}~,
\end{align}
and we are left to show that both terms \eqref{eq:sum1} and
\eqref{eq:sum2} converge to zero. 

The term
\eqref{eq:sum1} can be written as
$$
\eqref{eq:sum1} = \frac{1}{K_n} \sum_{k=1}^{K_n-1} \bigl [G_n^{-1} (t_{n,k}) -
G^{-1}(t_{n,k}) + a_{n,k} \bigr ]^2~,
$$
with $a_{n,k} \defeq G_n^{-1}(t_{n,k+1}) - G_n^{-1} (t_{n,k})$;
note that $\sum_{k=1}^{K_n-1} |a_{n,k}| \le C$.

Next, write
\begin{align}
 \frac{1}{K_n} \sum_{k=1}^{K_n-1} \bigl [G_n^{-1} (t_{n,k}) -
G^{-1}(t_{n,k}) + a_{n,k} \bigr ]^2
= & \; \; \; \; \;
\frac{1}{K_n} \sum_{k=1}^{K_n-1} \bigl [G_n^{-1} (t_{n,k}) -
G^{-1}(t_{n,k})\bigr ]^2 \label{eq:dec1} \\
& + \frac{2}{K_n} \sum_{k=1}^{K_n-1} \bigl [G_n^{-1} (t_{n,k}) -
G^{-1}(t_{n,k})\bigr ] \cdot a_{n,k} \label{eq:dec2} \\
& + 
\frac{1}{K_n} \sum_{k=1}^{K_n-1} a_{n,k}^2~. \label{eq:dec3}
\end{align}
The term on the right-hand side \eqref{eq:dec1} converges to zero by assumption. The
term \eqref{eq:dec3} converges to zero because
$$\sum_{k=1}^{K_n-1} a_{n,k}^2 \le \Bigl (\sum_{k=1}^{K_n-1} |a_{n,k}|
\Bigr )^2 \le C^2~.$$
Since both the term on the right-hand side of \eqref{eq:dec1} and the
term \eqref{eq:dec3} converge to
zero, the term \eqref{eq:dec2} converges to zero as well by the
Cauchy-Schwarz inequality. Consequently,
the term \eqref{eq:sum1} converges to zero.

By a completely analogous argument, the term \eqref{eq:sum2} converges to zero too.

We have thus established that the
convergence~\eqref{eq:conv-sum-finite} implies the convergence
\eqref{eq:conv-int}. By a similar argument, one can establish the
reverse fact that the
convergence~\eqref{eq:conv-int} implies the convergence
\eqref{eq:conv-sum-finite}.~\qed

\bigskip
It is useful to discuss Lemma~\ref{lem:grid-con} a bit further.
For two c.d.f.'s $G_1$ and $G_2$ on the real line, define
\begin{equation} \label{eq:p-norm}
||G_1 - G_2||_p \defeq \sqrt{\frac{1}{p}\sum_{i=1}^p \bigl [G_1^{-1}((i-0.5)/p) -
  G_2^{-1}((i-0.5)/p) \bigr ]^2}~.
\end{equation} 
Two results are noted. 

First, the left-hand expression of equation
\eqref{eq:conv-sum-finite} can be written as 
\begin{equation*} \label{eq:rewrite}
||G_n - G||_p^2
\end{equation*} 
when $p = K_n$. Since $p \to \infty$,
Lemma~\ref{lem:grid-con} states that, under the given set of assumptions,
\begin{equation} \label{eq:lem-alt}
\mbox{$G_n \Rightarrow G \quad$ if and only } \quad
||G_n - G||_p^2 \to 0~.
\end{equation} 

Second, a triangular inequality holds in the sense
that for three c.d.f.'s $G_1, G_2$, and $G_3$ on the real line,
\begin{equation}\label{eq:triangle}
||G_1 - G_2||_p\le ||G_1 - G_3||_p + ||G_2 - G_3||_p~.
\end{equation} 
This second fact follows since, for example, $\sqrt{p} \cdot ||G_1 - G_2||_p$ is
the Euclidian \mbox{distance} between the two vectors 
$(G_1^{-1}(0.5/p), \ldots, G_1^{-1}((p-0.5)/p))'$ and
$(G_2^{-1}(0.5/p), \ldots, \mbox{$G_2^{-1}((p-0.5)/p)$})'$.

These two results are summarized in the following corollary.

\begin{corollary} \label{cor:ap1}$\,$
\begin{itemize} 

\item[(i)]
Let $\{G_n\}$ and $G$ be c.d.f.'s on the real line and assume that there
exists a compact interval that contains the support of $G$ as well as the
support of $G_n$ for all $n$ large enough. For $0 < \alpha < 1$, let $G_n^{-1} (\alpha)$
denote an $\alpha$ quantile of $G_n$ and let $G^{-1}(\alpha)$ denote
an $\alpha$ quantile of $G$. Also assume that $p \to \infty$.

Then $G_n \Rightarrow G$ if and only if
$$
||G_n - G||_p^2 \to 0~,
$$
where $||\cdot||_p$ is defined as in \eqref{eq:p-norm}.

\item[(ii)] Let $G_1$, $G_2$, and $G_3$ be c.d.f.'s on the real
  line. Then
$$
||G_1 - G_2||_p\le ||G_1 - G_3||_p + ||G_2 - G_3||_p~.
$$
\end{itemize} 
\end{corollary}

\bigskip
{\sc Proof of Theorem~\ref{theo:individual}. }
As shown by \cite{silverstein:1995}, $F_n \Rightarrow F$ almost surely.
Therefore, by Corollary~\ref{cor:ap1}(i),
\begin{equation} \label{eq:ind1}
\frac{1}{p} \sum_{i=1}^p [\lambda_{n,i} - F^{-1}((i-0.5)/p)]^2 \toas 0~,
\end{equation} 
recalling that $\lambda_{n,i}$ is a $(i-0.5)/p$ quantile of $F_n$; see Remark~\ref{rem:quant}.
The additional fact that
\be \label{eq:ind2}
\frac{1}{p}\sum_{i=1}^p\left[q_{n,p}^i(\boldsymbol{\tau}_n) - F^{-1}((i-0.5)/p)\right]^2\stackrel{\rm
  a.s.}{\longrightarrow}0
\ee
follows from the \MP equation \eqref{eq:MP}, Lemma~A.2 of
\cite{ledoit:wolf:2012}, Assumption~(A3), the definition of
$q_{n,p}^i(\boldsymbol{\tau}_n)$, and Corollary~\ref{cor:ap1}(i) again. The convergences
\eqref{eq:ind1} and \eqref{eq:ind2} together with the triangular
inequality for the Euclidian distance in $\R^p$ then imply that
$$
\frac{1}{p} \sum_{i=1}^p \left [q_{n,p}^i(\boldsymbol{\tau}_n) -
  \lambda_{n,i} \right ]^2 \toas 0~,
$$
which is the statement to be proven.~\qed

\bigskip 
{\sc Proof of Theorem~\ref{theo:estimator}. }
For any probability measure $\widetilde H$ on the nonnegative real line
and for any $\wt c > 0$, let $F_{\widetilde H, \wt c}$ denote
the c.d.f.\ on the real line induced by the corresponding solution of the \MP
equation~\eqref{eq:MP}.
More specifically, for each $z \in \C^+$, $m_{F_{\widetilde H,\wt c}}(z)$ 
is the unique solution for $m \in \C^+$ to the equation
$$
m =
\int_{-\infty}^{+\infty}\frac{1}{\tau
\left[1-\wt c- \wt c\,z\,m \right]-z}\,d \widetilde H(\tau)~.
$$
In this notation,  $F = F_{H,c}$.

Recall that $F_n$ denotes the empirical c.d.f.\ of the sample
eigenvalues $\boldsymbol{\lambda}_n$. Furthermore, for $\tf \defeq
(t_1, \ldots, t_p)' \in [0,
\infty)^p$, denote by $\widetilde H_{\tf}$ the probability
distribution that places mass~$1/p$ at each of the $t_i \; (i = 1,
\ldots, p)$. 
The objective function in equation~\eqref{eq:optim} can then be
re-expressed as
$$
||F_{\widetilde H_{\tf}, \widehat c_n} - F_n||_p^2~,
$$
where $||\cdot||_p$ is defined as in \eqref{eq:p-norm}.
Note here that $F_{\widetilde H_{\tf}, \widehat c_n}$ is nothing
else than
$F_{n,p}^\mathbf{t}$ of equation~\eqref{eq:questF}; 
but for the purposes of this proof, the notation
$F_{\widetilde H_{\tf}, \widehat c_n}$ is more convenient.

Consider the following infeasible estimator of the
limiting spectral distribution~$H$:
\begin{equation} \label{eq:H-infeasible}
\overline H_n \defeq \argmin_{\widetilde H} ||F_{\widetilde H,
  \widehat c_n} - F_n||_p^2~,
\end{equation} 
where the minimization is over {\it all} probability measures
$\widetilde H$ on the real line; the estimator $\overline H_n$ is
infeasible, since one cannot minimize over all probability measures on
the real line in practice.
By definition,
\begin{equation} \label{eq:trivial}
||F_{\overline H_n, \wh c_n}-F_n||_p \le 
||F_{H, \wh c_n}-F_n||_p~.
\end{equation} 
Therefore,
\begin{align*}
  ||F_{\overline H_n, \wh c_n} - F||_p & \le 
  ||F_{\overline H_n, \wh c_n} - F_n||_{p} + ||F_n - F||_{p} 
  \quad \mbox{(by Corollary~\ref{cor:ap1}(ii))} \\
& \le 
  ||F_{H, \wh c_n} - F_n||_{p} + ||F_n - F||_{p}\quad  \mbox{ (by
  \eqref{eq:trivial}} \\
  & \le  
  ||F_{H, \wh c_n} - F_{H, _c}||_{p} + ||F_{H, c} - F_n||_{p} + ||F_n -
  F||_{p} \quad \mbox{ (by Corollary~\ref{cor:ap1}(ii))}\\
  & =  
  ||F_{H, \wh c_n} - F||_{p} + 2\, ||F_n - F||_{p} \quad \mbox{ (since
    $F_{H,c} = F$)} \\
  & \eqdef A + B~.
\end{align*}
In the case $c < 1$, combining Corollary~\ref{cor:ap1}(i) 
with Lemma A.2 of
\cite{ledoit:wolf:2012} shows that $A \to 0$ almost surely.
In the case $c > 1$, one can also show that $A \to 0$ almost surely: 
Lemma A.2 of \cite{ledoit:wolf:2012} implies that 
$\underline F_{H,\widehat c_n} \Rightarrow \underline F$ almost surely; then use
equation~\eqref{eq:mF2} together with the fact that $\widehat c_n \to c$
to deduce that also $F_{H,\widehat c_n} \Rightarrow F$ almost surely; finally
apply Corollary~\ref{cor:ap1}(i).
Combining Corollary~\ref{cor:ap1}(i) with the fact that $F_n \Rightarrow
F$ almost surely~(\citeauthor{silverstein:1995}, 1995) shows that $B \to 0$
almost surely\ in addition to $A \to 0$ almost surely. Therefore, $||F_{\overline H_n,
  \wh c_n} - F||_p  \to 0$ almost surely.
Using Corollary~\ref{cor:ap1}(i) again shows that 
$F_{\overline H_n, \wh c_n} \Rightarrow F$~almost surely.

A feasible estimator of $H$ is given by 
\begin{equation*} \label{eq:H-feasible}
\widehat H_n \defeq \argmin_{\widetilde H_{\tf} \in {\mathcal P}_n} ||F_{\widetilde H_{\tf},
  \widehat c_n}-F_n||_p^2
\end{equation*} 
instead of by \eqref{eq:H-infeasible}, where the subset $\mathcal{P}_n$
denotes the set of probability measures that are
equal-weighted mixtures of $p$ point masses on the nonnegative real line:
\begin{equation*}\label{eq:P-n}
{\mathcal P}_n \defeq \Bigl \{\widetilde H_{\tf}: 
\widetilde{H}_{\tf}(x) \defeq \frac{1}{p} \sum_{i=1}^p\bone_{\{x\geq t_i\}}~, \mbox{
  where } \tf \defeq (t_1, \ldots, t_p)^\prime \in [0, \infty)^p \Bigr
\}~.
\end{equation*} 
The fact that the minimization
over a finite but dense family of probability measures, instead of
all probability measures on the nonnegative real line, does not affect the strong
consistency of the estimator of~$F$ follows by arguments similar to those used in the proof of
Corollary~5.1(i) of \cite{ledoit:wolf:2012}. Therefore, it also holds that
$F_{\widehat H_n, \wh c_n} \Rightarrow F$ almost surely.

\smallskip

Having established that $F_{\widehat H_n, \wh c_n} \Rightarrow F$
almost surely, it follows that also $\widehat H_n \Rightarrow H$
almost surely; see the proof of
Theorem~5.1(ii) of \cite{ledoit:wolf:2012}. Since $\widehat H_n$ is
recognized as
the empirical distribution (function) of the $\widehat \tau_{n,i} \;
(i = 1, \ldots, p)$, $\widehat \tau_{n,i}$ is
a $(i-0.5)/p$ quantile of~$\widehat H_n$; see Remark~\ref{rem:quant}.
Therefore, it follows from Corollary~\ref{cor:ap1}(i) that
\begin{equation} \label{eq:tau1}
\frac{1}{p} \sum_{i=1}^p [\widehat \tau_{n,i} - H^{-1}((i-0.5)/p)]^2 \toas 0~,
\end{equation} 
The additional fact that
\be \label{eq:tau2}
\frac{1}{p}\sum_{i=1}^p\left[\tau_{n,i} - H^{-1}((i-0.5)/p)\right]^2\stackrel{\rm
  a.s.}{\longrightarrow}0
\ee
follows directly from Assumption (A3) and Corollary~\ref{cor:ap1}(i) again. The convergences
\eqref{eq:tau1} and \eqref{eq:tau2} together with the triangular
inequality for the Euclidian distance in $\R^p$ then imply that
$$
\frac{1}{p} \sum_{i=1}^p \left [\widehat \tau_{n,i} - \tau_{n,i}\right ]^2 \toas 0~,
$$
which is the statement to be proven.~\qed

\bigskip
{\sc Proof of Theorem~\ref{prop:oracle-cons}. }
The claim for the case $p < n$
follows immediately from Proposition~4.3(ii) of
\cite{ledoit:wolf:2012}. 

\smallskip

To treat the case $p >n$, let $j$ denote the smallest integer for which
$\lambda_i > 0$. Note that $(j-1)/p \to (c-1)/c$ almost surely by the results
of \cite{bai:silverstein:1999}; indeed, since 
the~$\lambda_i$ are sorted in increasing order, $(j-1)/p$
is just the fraction of sample eigenvalues that are equal to zero.

Now restrict attention to the set of probability one
on which $\widehat{\breve m_{\underline F} (0)} \to 
\breve m_{\underline F} (0)$, $\widehat H_n \Rightarrow H$,
and $(j-1)/p \to (c-1)/c$.
Adapting Proposition~4.3(i)(a) of \cite{ledoit:wolf:2012} to the
continuous part of $F$, it can be shown that $\breve m_{\wh H_n, \wh c_n} (\lambda) \to \breve
m_F(\lambda)$ uniformly in $\lambda \in \supp(\underline F)$, except for two
arbitrarily small regions at the lower and upper end of
$\supp(\underline F)$. We can write
\begin{align}
||\widehat S_n - S_n^{or}||_F^2  & =  
\frac{j-1}{p}
\biggl ( \frac{1/\wh c_n}{(1-1/\wh c_n) \, \wh{\breve m_{\underbar{\mbox{\tiny $F$}}} (0)}}
- \frac{1/c}{(1-1/c) \, {\breve m_{\underbar{\mbox{\tiny $F$}}}
    (0)}} \bigg )^2 \nonumber \\
& + \frac{1}{p} \sum_{i=j}^p \biggl (
\frac{\lambda_i}{\bigl | 1 - \widehat c_n - \widehat c_n \, \lambda_i \,
\breve m_{F_{\wh H_n, \wh c_n}}(\lambda_i) \bigr |^2} -
\frac{\lambda_i}{\bigl | 1 - c - c \, \lambda_i \,
\breve m_F(\lambda_i) \bigr |^2} \biggr )^2 \nonumber \\
& \eqdef D_1 + D_2~. \label{eq:d1-d2}
\end{align}
The fact that the summand~$D_1$ converges to zero is obvious,
keeping in mind that $c > 1$ and $\breve
m_{\underbar{\mbox{\tiny $F$}}} (0) > 0$.
The fact that the summand $D_2$ converges to zero follows by arguments
similar to those in the proof of Proposition~4.3(i)(b) of
\cite{ledoit:wolf:2012}. 

We have thus shown that there exists a set of probability one on which
$||\widehat S_n - S_n^{or}||_F \to 0$.~\qed

\section{Justification of Remark \protect\ref{rem:total-var}}
\label{app:offdiagonal}

{\sc Notation. }

\begin{packed_itemize}
\item Let
$y$ be a real $p$-dimensional random vector with covariance
matrix $\Sigma$. 
\item
Let $I_k$ denote the $k$-dimensional identity matrix, where $1\leq k\leq p$. 
\item
Let $W$ be a real nonrandom matrix of dimension $p\times k$ such that $W'W=I_k$. 
\item
Let $w_i$ denote the $i$th column vector of $W$ ($i=1,\ldots,k$).
\end{packed_itemize}
We start from the following two statements.
\begin{enumerate}
\renewcommand{\labelenumi}{(\arabic{enumi})}
\item If $\Cov[w_i'y,w_j'y]=0$ for all $i\neq j$ then the variation attributable to the set of random variables $(w_1'y,\ldots,w_k'y)$ is $\sum_{i=1}^k\Var[w_i'y]$.\label{stat:diagonal}
\item If $R$ is a $k\times k$ rotation matrix, that is,~$R'R=RR'=I_k$,
and $\widetilde{w}_i$ is the $i$th column vector of the matrix $WR$, then the variation attributable to the rotated variables $(\widetilde{w}_1'y,\ldots,\widetilde{w}_k'y)$ is the same as the variation attributable to the original variables $(w_1'y,\ldots,w_k'y)$.\label{stat:rotation}
\end{enumerate}
Together, Statements (\ref{stat:diagonal}) and (\ref{stat:rotation})
imply that, even if $\Cov[w_i'y,w_j'y]\neq0$ for $i\neq j$, the variation
attributable to $(w_1'y,\ldots,w_k'y)$ is still $\sum_{i=1}^k\Var[w_i'y]$.

\medskip
{\sc Proof. }
Let us choose $R$ as a matrix of eigenvectors of
$W'\Sigma W$. Then $(WR)'\Sigma(WR)$ is diagonal and
$(WR)'(WR)=I_k$. Therefore, by Statement (\ref{stat:diagonal}), the
variation attributable to $(\widetilde{w}_1'y,\ldots,\widetilde{w}_k'y)$ is $\sum_{i=1}^k\Var[\widetilde{w}_i'y]=\tr[(WR')\Sigma(WR)]$. By the properties of the trace operator, this is equal to $\tr(W'\Sigma
W)=\sum_{i=1}^k\Var[w_i'y]$. By Statement (\ref{stat:rotation}), it is the same as the
variation attributable to $(w_1'y,\ldots,w_k'y)$.~\qed
\end{appendix}

\end{document}